\documentclass[12pt]{article}
\usepackage{amsmath}
\newtheorem{thm}{Theorem}[section]
\newtheorem{lem}[thm]{Lemma}
\newtheorem{cor}[thm]{Corollary}
\newtheorem{prop}[thm]{Proposition}
\newenvironment{proof}{\noindent{\it Proof.}}{$\hfill \Box$}

\newtheorem{defn}[thm]{Definition}
\newtheorem{exa}[thm]{Example}
\newtheorem{rem}[thm]{Remark}
\usepackage{amsfonts}
\usepackage{amssymb}
\usepackage{graphicx}
\usepackage{float}
\usepackage{graphicx}
\usepackage{caption}
\usepackage{subcaption}
\usepackage{tikz}
\usetikzlibrary{calc} 
\tikzstyle{vertex}=[circle, fill=black, inner sep=.3pt, draw]

\usepackage{graphicx}
\usepackage{caption}
\usepackage{subcaption}
\usetikzlibrary{backgrounds}
\usetikzlibrary{arrows}
\usetikzlibrary{shapes,shapes.geometric,shapes.misc}

\tikzstyle{tikzfig}=[baseline=-0.25em,scale=0.5]

\pgfkeys{/tikz/tikzit fill/.initial=0}
\pgfkeys{/tikz/tikzit draw/.initial=0}
\pgfkeys{/tikz/tikzit shape/.initial=0}
\pgfkeys{/tikz/tikzit category/.initial=0}

\pgfdeclarelayer{edgelayer}
\pgfdeclarelayer{nodelayer}
\pgfsetlayers{background,edgelayer,nodelayer,main}

\newcommand{\tikzfig}[1]{%
{\tikzstyle{every picture}=[tikzfig]
\IfFileExists{#1.tikz}
  {\input{#1.tikz}}
  {%
    \IfFileExists{./figures/#1.tikz}
      {\input{./figures/#1.tikz}}
      {\tikz[baseline=-0.5em]{\node[draw=red,font=\color{red},fill=red!10!white] {\textit{#1}};}}%
  }}%
}

\usepackage{times}

\usepackage{times}
\usepackage{rotating}
\usepackage[T1]{fontenc}

\begin{document}
	\pagenumbering{gobble}
	
	\Large
	\begin{center}
	 Cayley graphs and $G$-graphs of Gyro-groups	
		\\
	
		\hspace{10pt}

		\large

		{Neda Moradi}$^{a}$, {Gholam Hossein Fath-Tabar}$^{a,*}$,{ Alain Bretto}$^{b}$ 
		
		\hspace{10pt}

	\small

		$^{a}${Department of Pure Mathematics, Faculty of Mathematical Sciences, University of Kashan, Kashan, I. R. Iran}\\
		
     $^{b}${ Universite de Caen, GREYC CNRS UMR-6072, Campus II Bd Marechal Juin BP 5186, 14032 Caen cedex Caen, France}

	\end{center}
	
	\hspace{10pt}
	
	\normalsize
		
\begin{center}
\textbf{This work is in memory of Professor Ali Reza Ashrafi, an outstanding teacher of science and ethics.}
\end{center}

\begin{abstract}
In this paper the structure of the Cayley graphs and $G$-graphs of some gyro-groups are studied and some properties of them will be proved.  Moreover we review some special gyro-groups including: gyro-commutative gyro-groups, dihedral gyro-groups and dihedralized gyro-groups. Then we try to establish some important properties of their associated $G$-graphs. Finally we will examine the symmetry of Cayley graphs and $G$-graphs of some gyro-groups.

\begin{flushleft}

\end{flushleft}		

\end{abstract}

\section{Preliminaries}
The concept of \textbf{gyro-groups } was first introduced by Ungar as an extension of the study of Lorentz groups.
In 1988 in \cite{Ungar}. Ungar presented the parametric expression of Lorentz transformations in terms of acceptable
 relative velocities and directions. The group structure resulting from the realization of the Lorentz group,
 together with Einstein's speed addition law, enabled him to establish a profound structural connection between Thomas's change of direction and Einstein's addition law. 

Ungar extended Thomas's change of direction and referred to its generalized form as Thomas rotation. He showed that rotations are  automorphisms and exhibit group-like properties, which he
 called the gyro-group. Gyro-groups represent a particular generalization of groups that do not possess
 associativity. 
 
 We start the discussion with the simplest algebraic structure which is the groupoid. Indeed a groupoid is a set having one binary operation satisfying only closure.  
 
\begin{defn}
A groupoid $ (G, \oplus)$ is called a gyro-group if its binary operation
satisfies the following axioms:
\begin{itemize}
\item[(G1)]there is $0 \in G $such that $0 \oplus a = a$ for all $a \in G$;
\item[(G2)] for any $ a \in G$, there is $b \in G$ such that $b \oplus a = 0$;
\item[(G3)] for any $a, b \in G$, there is an automorphism $gyr[a, b]: G \longmapsto G$ such that for
any $c \in G$,
$a \oplus (b \oplus c) = (a \oplus b)  \oplus gyr[a, b]c$;
\item[(G4) ]for any $a, b \in G$, $gyr[a, b] = gyr[a \oplus b, b]$.
\end{itemize}
\end{defn}

It is recommended to refer to \cite{Sabinin, Suksumran2015, Suksumran2016, Maungchang2021} gyro-groups for more details.

Cayley graphs are a good tool for group representation. They are always regular, but they do not give us much information about the corresponding group, in addition, two  isomorphic groups have isomorphic Cayley graphs, but converse is not true. But $G$-graphs as a new representation of groups do not have these limitations. 
As we know that the group is also a gyro-group, we will try to study the structure of Cayley graph and G-graph for gyro-groups. The concept of Cayley graph associated with the gyro-groups was first investigated in \cite{Bussaban} by Bussaban and et, al and a number of the characteristics of Cayley graphs of gyro-groups were discussed. The structure of $G$-graph associated to a finite group $G$ was first introduced by Alain Bretto, et. al., in 2005 in \cite{Beretto2005}. This graph is very similar to Cayley graph in some ways, while not having some of their limitations.
It is well-known that the $G$-graph of a groups is connected if and only if the generating set spans a group. However, this fact need not be satisfied for the $G$-graph of  gyro-groups. We show that the spanning condition does not guarantee connectedness of $G$-graphs of gyro-groups.

 In this article we are going to construct the structure of the $G$-graphs of some finite gyro-groups. Moreover,  we continue to prove some properties of them, including the connectivity.   Finally we check the symmetry of $G$-graphs of gyro-groups and we discus about the conditions of symmetry  $G$-graphs of gyro-groups.

\begin{defn} \cite{Bretto2007} 
 Let $ (G, S)$ be a group with a set of generators
$$S = \lbrace s_{1}, s_{2}, s_{3}\cdots s_{k}\rbrace, k \geq 1. $$
 For any $ s \in S$, we consider the left action of the subgroup $H= \langle s \rangle $ on $G$. So we have a partition $G=\sqcup  \langle s \rangle x $ , $ x \in T_{s}$, where $T_{s}$ is a right  transversal of $ \langle s \rangle $. The cardinality of $  \langle s \rangle $ is $o(s)$, the order of the element s. Let us consider the cycles
$$(s)x = (x, sx, s^{2}x \cdots s^{o(s)-1}x)$$
of the permutation $g_{s}: x \longmapsto sx $ of $G$.  We define a graph denoted by

$\Phi(G; S) = (V; E; \epsilon )$ in the following way:

\begin{itemize}
\item[I)] The vertices of $ \Phi(G; S)$ are the cycles of $g_{s}$, $s \in S, $ i.e.,
  $V = \sqcup V_{s}$ , $s \in S$ with $V_{s} = \lbrace (s)x, x \in T_{s} \rbrace$.
\item[II)] For two different vertices $(s)x, (t)y \in V$, that $card( \langle s \rangle x \bigcap \langle t \rangle y) = p$,

 $ p \geq 1$, then $ {(\langle s \rangle x , \langle t \rangle y)}$ is a $p$-edge. 
\end{itemize}

$\Phi(G; S)$ is called graph from groups or $G$-graph. 
\end{defn}
Interesting features of G-graphs are described in detail in references

 \cite{Bretto2011,Bretto2004,Bretto-2007,Beieke,Ashrafi2019}.

\subsection{$G$-graph of gyro-groups}

\begin{defn}
	Let $ (G, S)$ be a gyro-group with a set of generators
	$S = \lbrace s_{1}, s_{2}, s_{3}\cdots s_{k}\rbrace$, $k \geq 1$.
	For any $ s \in S$, we consider the left action of the sub-gyro-group $H= \langle s \rangle $ on $G$. So we have a partition $G=\sqcup  \langle s \rangle \oplus x $ , $ x \in T_{s}$, where $T_{s}$ is a right  transversal of $ \langle s \rangle $. The cardinality of $  \langle s \rangle $ is $o(s)$, the order of the element s. Let us consider the cycles

	\begin{flushleft}
		$(s) \oplus x = (x, s \oplus x, (s\oplus s)\oplus x, \cdots \oplus ,
	\underbrace { (s \oplus \cdots \oplus s)}_{o(s)-1 ~~~ times} \oplus x ) $  		
	\end{flushleft}
	
	of the permutation $g_{s}: x \longmapsto s \oplus x $ of $G$.  We define a graph denoted by 
	$\Phi(G; S) = (V; E; \epsilon )$ in the following way:
	
	\begin{itemize}
		\item[I)] The vertices of $ \Phi(G; S)$ are the cycles of $g_{s}$, $s \in S, $ i.e.,
		$V = \sqcup V_{s}$ , $s \in S$ with $V_{s} = \lbrace (s) \oplus x, x \in T_{s} \rbrace$.
		\item[II)] For two different vertices $(s) \oplus x, (t) \oplus y \in V$, that 
		
		$card( \langle s \rangle \oplus x \bigcap \langle t \rangle \oplus y) = p$, $ p \geq 1$, then $ {(\langle s \rangle \oplus x , \langle t \rangle \oplus y)}$ is a $p$-edge. 
	\end{itemize}
	
$\Phi(G;S)$ is called
graphs from gyo-groups or\textbf{ $G$-gyro-graph}. 	
\end{defn}

Here there is an example of the structure of the  $G$-graph of  the gyro-group $G_{8}$, with the Cayley table, gyro table and figure of the graph with respect to the generating set $S=\{1, 2\}$. 

\begin{exa} \cite{Farzaneh2022}\label{ma}
Consider the gyro group $G_{8}= \lbrace 0,1,2,3,4,5,6,7 \rbrace $ and

 $S=\lbrace 1,2 \rbrace$ with $  A = (1 6)(2 5) $, then the gyro table of it is as follows. 
\begin{table}
\begin{center}
\small{
\begin{tabular}{c|llllllll lc|llllllll}
$\oplus$ &  $0$ & $1$ & $2$ & $3$ & $4$ & $5$ & $6$ & $7$ & & gyro &  $0$ & $1$ & $2$ & $3$ & $4$ & $5$ & $6$ & $7$ \\
\hline 
$0$ & $0$ & $1$ & $2$ & $3$ & $4$ & $5$ & $6$ & $7$ & & $0$ & $I$ & $I$ & $I$ & $I$ & $I$ & $I$ & $I$ & $I$\\
$1$ & $1$ & $0$ & $3$ & $2$ & $5$ & $4$ & $7$ & $6$ &  & $1$ & $I$ & $I$ & $A$ & $A$ & $A$ & $A$ & $I$ & $I$\\
$2$ &  $2$ & $3$ & $ 0$ & $1$ & $6$ & $7$ & $4$ & $5$ & & $2$ &  $I$ & $A$ & $I$ & $A$ & $A$ & $I$ & $A$ & $I$ \\
$3$ &  $3$ & $5$ & $6$ & $0$ & $7$ &  $1$ & $2$ & $4$  & & $3$ &  $I$ & $A$ & $A$ & $I$ & $I$ & $A$ & $A$ & $I$ \\
$4$ &  $4$ & $2$ & $1$ & $7$ & $0$ & $6$ & $5$ & $3$  &  & $4$ &  $I$ & $A$ & $A$ & $I$ & $I$ & $A$ & $A$ & $I$ \\
$5$ &  $5$ & $4$ & $7$ & $6$ & $1$ & $0$ & $3$ & $2$  &  & $5$ &  $I$ & $A$ & $I$ & $A$ & $A$ & $I$ & $A$ & $I$ \\
$6$ &  $6$ & $7$ & $4$ & $5$ & $2$ & $3$ & $0$ & $1$  &  & $6$ &  $I$ & $I$ & $A$ & $A$ & $A$ & $A$ & $I$ & $I$ \\
$7$ &  $7$ & $6$ & $5$ & $4$ & $3$ & $2$ & $1$ & $0$   & & $7$ &  $I$ & $I$ & $I$ & $I$ & $I$ & $I$ & $I$ & $I$ 
\end{tabular}
}
\caption{ The gyro table of $G_8$}
\end{center}
\end{table}

Obviously $ S $ is a generating set for $ G_{8} $,  because:
\begin{flushleft}
$ 3= 1\oplus 2 $, $4=1 \oplus ( (1 \oplus 2) \oplus 1)  $, $ 5=(1\oplus 2) \oplus 1 $, $ 6=(1 \oplus2) \oplus 2 $ and $ 7=2 \oplus ( (1\oplus 2) \oplus 1 $.
\end{flushleft}
 We compute the $2$-cycles of $V_{r} = \lbrace (r) \oplus x; x \in G_{8}, r \in S \rbrace$. Then 

\begin{flushleft}
$ V_{1} = \lbrace (1) \oplus x= (x, 1 \oplus x); x \in G_{8} \rbrace= \lbrace  (0, 1), (2, 3), (4, 5), (6, 7)\rbrace $.
\end{flushleft}

Similarly for  $V_{2}  $ we have:

\begin{flushleft}
$ V_{2} = \lbrace (2) \oplus y= (y, 2 \oplus y); y \in G_{8} \rbrace=\lbrace (0, 2), (1, 3), (4, 6), (5, 7) \rbrace $.
\end{flushleft}

We know that $ V(\Phi(G_{8}; S = \lbrace1; 2 \rbrace)) = V_{1} \bigcup V_{2} $, as follows:

\begin{center}
$ V_{1} \bigcup V_{2}= \lbrace  (0, 1), (2, 3), (4, 5), (6, 7)\rbrace \bigcup  \lbrace (0, 2), (1, 3), (4, 6), (5, 7) \rbrace  $
\end{center}

\begin{figure}[H]
\begin{center}
\begin{tikzpicture}
\node (v1) at (0,0) {$(0,1)$};
\node (v2) at (2,0) {$(2,3)$};
\node (v3) at (2,-2) {$(1,3)$};
\node (v4) at (0,-2) {$(0,2)$};
\draw (v1) edge (v3);
\draw (v1) edge (v4);
\draw (v2) edge (v3);
\draw (v2) edge (v4);
\end{tikzpicture}
\hspace{.5cm}
\begin{tikzpicture}
\node (v1) at (0,0) {$(4,5)$};
\node (v2) at (2,0) {$(6,7)$};
\node (v3) at (2,-2) {$(5,7)$};
\node (v4) at (0,-2) {$(4,6)$};
\draw (v1) edge (v3);
\draw (v1) edge (v4);
\draw (v2) edge (v3);
\draw (v2) edge (v4);
\end{tikzpicture}
\caption{ $\Phi(G_8 , S=\{1, 2 \})$}
\end{center}
\end{figure}
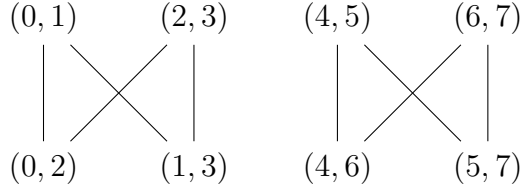

 As you see in Figure 1, the $G$-graph of this gyro-group is disconnected with two isomorphic connected components of $C_4$. Note that although $ S $ is a generating set for $ G $,  $ \Phi(G_{8}, S) $  is not connected. Just like the case for Cayley graph of gyro-groups.

\end{exa}

 \begin{defn}\cite{Bussaban}
  Let $ S $ be a subset of a gyro-group $ G $. The set $ S $ is said to be symmetric if for each element $ s \in S $, $ \ominus s \in S $. The left-generating set by $ S $, written by $ ( S \rangle$, is
$$( S \rangle = \{ s_{n} \oplus (\cdots \oplus ( \cdots s_{3} \oplus (s_{2} \oplus s_{1})) \cdots )\vert s_{1}, s_{2}, \cdots s_{n}  \in S \} .$$
If  $ (S \rangle = G $, we say that $ S $ left-generates  $ G $, or $ G $ is left-generated by $ S $. The right-generating set is defined in a similar fashion \cite{Maungchang}.
 \end{defn}

\begin{exa}
Consider the gyro-group $ G_{8}$ with the generating set  $S = \lbrace 1, 3 \rbrace $.  According to the explanation of the previous example, we have:
\begin{flushleft}
$ V(\Phi(G_{8}, S = \lbrace1, 3 \rbrace )) = V_{1} \bigcup V_{3}$

$= \lbrace (1) \oplus x = (x, 1 \oplus x) \rbrace \vert x \in G_{8}\rbrace  \bigcup  \lbrace (3) \oplus y = (y, 3 \oplus y) \rbrace \vert y \in G_{8} \rbrace$
$= \lbrace (0, 1), (2, 3), (4, 5), (6, 7)\rbrace \bigcup  \lbrace (0, 3), (1, 5), (2, 6), (4, 7) \rbrace $.
\end{flushleft}
We can see that $ \Phi (G_{8}, S) $  is a connected $G$-graph of $G_8$ isomorphic to a cycle $ C_{8} $.  

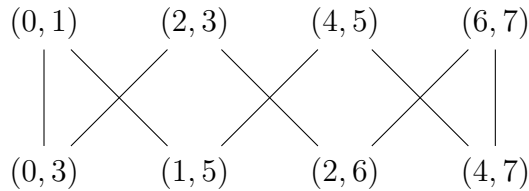
\begin{figure}[H]
 \begin{center}
\begin{tikzpicture}
\node (v1) at (0,0) {$(0,1)$}; \node (v2) at (2,0) {$(2,3)$}; \node (v3) at (4,0) {$(4,5)$}; \node (v4) at (6,0) {$(6,7)$}; 
\node (v5) at (0, -2) {$(0,3)$}; \node (v6) at (2,-2) {$(1,5)$}; \node (v7) at (4, -2) {$(2,6)$}; \node (v8) at ( 6,-2) {$(4,7)$}; 
\draw (v1) -- (v5); \draw (v1) -- (v6); \draw (v2) -- (v5); \draw (v2) -- (v7); \draw (v3) -- (v6); \draw (v3) -- (v8); \draw (v4) -- (v7); 
\draw (v4) -- (v8);
\end{tikzpicture}
\caption{$\Phi(G_8, S= \{ 1, 3\})$.}
\end{center}
\end{figure}

Note that here $G_{8} = (S \rangle $.  In other words, $ S $ is a left generating set for $G_{8}$. (See \cite{Bussaban, Maungchang} for more details on the left and right generator sets).
\end{exa}

Therefore, we can state the following theorem.
\begin{thm} \cite{Farzaneh2022}
Suppose $ G $ is a gyro-group and $ S $ is a non-empty subset of $G$. Then the $G$-graph  $\Phi (G, S)  $ is connected if and only if $(S\rangle =G  $.
\end{thm}

\begin{proof}
First, let $ S $ be the left generating set of $G$ of size one. We show that for every $ x \in G $, there is a path from $ e $ to  $ x $ and from $ x $ to $ e $. Since $ (S\rangle= G $, there exist $s_{1}, s_{2}, s_{3},\cdots, s_{n} \in S$ such that 
$$x= s_{n}\oplus (\cdots \oplus  (s_{3}\oplus(s_{2}\oplus s_{1}))).$$
 But we have assumed that $ S $ has only one element, which means
 
  $x=s =s\oplus e= e \oplus s  $. Thus the first part is proved. 

Now suppose $ \Phi (G, S) $ is connected. It is obvious that $ (S\rangle \subseteq G $. Because of connectivity of $ \Phi (G, S)$,  
 there exists $ s \in S $ such that $ y= s \oplus e $, then $ y \in (s\rangle $ and $ G \subseteq (S\rangle $, which implies $ (S\rangle = G $.
 
 Now assume that  $card(S) \geq 2 $. For $s, s^{\prime} \in S$, we show that there is a path from the set of vertices $V_s$ to the set of vertices of $V_{s^{\prime}}$.  For $x, y \in G$ consider two vertices  $(s) \oplus x \in V_{s}$ and
$(s^{\prime}) \oplus y \in V_{s^{\prime}}$. Since $G = (S\rangle $, there are  $s_{1}, s_{2}, s_{3},\cdots, s_{n} \in S $ such that 
$$y =  (s_{n}\oplus (\cdots \oplus  (s_{3}\oplus(s_{2}\oplus s_{1}))\cdots ) ) \oplus x. $$
Moreover we have that
{\small{
\begin{eqnarray*} 
x  & \in & \left\langle s \right\rangle  \oplus x  \bigcap   (s_{1} \oplus x),\\
s_{1} \oplus x  & \in &  \left\langle s_{1} \right\rangle  \oplus x \bigcap  \left\langle  s_{0} \right\rangle  \oplus ( s_{1} \oplus x ),\\
(s_2 \oplus  s_1 )  \oplus x  & \in & \left\langle  ( s_2 \oplus  s_1 )  \right\rangle  \oplus x \bigcap \left\langle  s_{1} \right\rangle  \oplus  ((s_2 \oplus  s_1 )  \oplus x )),\\
(s_{n-1} \oplus \cdots \oplus  (s_{3} \oplus (s_{2}  \oplus s_{1}))) \oplus x) & \in & \left\langle   (s_{n-1} \oplus \cdots \oplus ( s_{3} \oplus (s_{2}  \oplus s_{1})) \oplus x )  \right\rangle  \\ 
& \bigcap & \left\langle   s_{n} \right\rangle \oplus  (s_{n-1} \oplus \cdots \oplus ( s_{3} \oplus (s_{2}  \oplus s_{1})) \oplus x ), \\
y & \in &  \left\langle  (s_{n-1} \oplus \cdots \oplus ( s_{3} \oplus (s_{2}  \oplus s_{1})) \oplus x \right\rangle \bigcap  \left\langle  s^{\prime} \right\rangle  \oplus y. 
\end{eqnarray*}	}}
	Therefore we could find  a path from $(s) \oplus x  \in  V_{s}$ to $(s^{\prime}) \oplus y  \in V_{s^{\prime}}$. So $ \Phi(G, S)$ is a connected graph.
	
	Conversely, assume $ \Phi(G, S) $ is connected and $ x \in G $. There exists $s_{i_{1}} \in S$ and $x_{i_{1}} \in T_{s_{i_{1}}}$, such that $x \in (s_{i_{1}}) \oplus x_{i_{1}} $; 
	hence $x =\underbrace{(s_{i_{1}} \oplus \cdots \oplus s_{i_{1}})}_{t_{i_{1}} times} \oplus   x_{i_{1}} $

Therefore we could find  a path from $(s) \oplus x  \in  V_{s}$ to $(s^{\prime}) \oplus y  \in V_{s^{\prime}}$. So $ \Phi(G, S)$ is a connected graph.

Conversely, assume $ \Phi(G, S) $ is connected and $ x \in G $. There exists $s_{i_{1}} \in S$ and $x_{i_{1}} \in T_{s_{i_{1}}}$, such that $x \in (s_{i_{1}}) \oplus x_{i_{1}} $; 
 hence $x =\underbrace{(s_{i_{1}} \oplus \cdots \oplus s_{i_{1}})}_{t_{i_{1}} times} \oplus   x_{i_{1}} $. 
 
 We know that  $e \in  (s_{i_{1}})$ and the graph is connected if there exists a path from $ x \in (s_{i_{1}}) \oplus x_{i_{1}}  $ to $e \in ((s_{i_{1}})$; hence

\begin{eqnarray*}	
 x & = & \underbrace{(s_{i_{1}} \oplus \cdots \oplus s_{i_{1}})}_{t_{i_{1}} times}  \oplus  x_{i_{1}}  \\
   x_{i_{1}} & = & \underbrace{(s_{i_{2}} \oplus \cdots \oplus s_{i_{2}})}_{t_{i_{2}} times}  \oplus  x_{i_{2}}, \cdots , x_{i_{k-2}}= \underbrace{(s_{i_{k-1}} \oplus \cdots \oplus s_{i_{k-1}})}_{t_{i_{k-1}} times}  \oplus  x_{i_{k-1}} \\
  x_{i_{k-1}} & =  &\underbrace{(s_{i_{k}} \oplus \cdots \oplus s_{i_{k}})}_{t_{i_{k}} times}  \oplus  x_{i_{k}}  
\end{eqnarray*}	 

But $x_{i_{k}} = e$, so $x =\underbrace{(s_{i_{1}} \oplus \cdots \oplus s_{i_{1}})}_{t_{i_{1}} times} \oplus \underbrace{(s_{i_{2}} \oplus \cdots \oplus s_{i_{2}})}_{t_{i_{2}} times}\oplus \cdots \oplus \underbrace{(s_{i_{k}} \oplus \cdots \oplus s_{i_{k}})}_{t_{i_{k}} times} $

which shows that every element $x$ of the gyro-group $G$ can be written by the elements of the set $S$, this implies  that $G$ is left generated by $S$.  
\end{proof}

\section{Main Result}
In this section we are going to bring the structure of some gyro-groups and their associated Cayley graph and $G$-graph of them. First of all we focus on the characterization of the gyro-groups of order $8$ distinguished by Ashrafi., et, al, in \cite{Ashrafi}, \cite{Ashrafi2010}, \cite{Ashrafi2022}.  Furthermore we try to obtain some properties of their $G$-graphs and Cayley graphs.

\subsection{Cayley graph and $G$-graph of gyro-groups of order 8}
There are exactly six gyro-groups of order $8$, which are presented in \cite{Ashrafi2022}, whose gyro-tables are in reference \cite{Mahdavi}.
Note that the gyro-groups mentioned in  \cite{Ashrafi2010}, i.e
. $K(1)$, $L(1)$, $M(1)$, $N(1)$, $O(1)$  are isomorphic with gyro-groups of order $8$ in paper \cite{Ashrafi2022}. We can see that \\
 $ G_{8, 1} \cong L(1) $, $ G_{8, 2} \cong K(1) $,  $ G_{8, 3} \cong O(1) $, $ G_{8, 5} \cong M(1) $, and 
$ G_{8, 6} \cong N(1) $. Among them the structure of the gyro-groups   of  $ G_{8, 1}$,  $ G_{8, 3}$ and  $ G_{8, 5}$ are gyro-commutative.

We have drawn all the $G$-graphs of the above gyro-groups and by analyzing the corresponding $G$-graphs, we have reached interesting results that we will mention below.
\begin{table}[H]\begin{center}
\begin{tabular}{c|cccccccc}
${\oplus}_{K(1)}$ & $0$ & $1$ & $2$ & $4$ & $5$ & $6$ & $7$ & $8$\\
\hline 
$0$ & $0$ & $1$ & $2$ & $2$ & $3$ & $4$ & $5$ & $6$ \\
$1$ & $1$ & $2$ & $3$ & $4$ & $5$ & $6$ & $7$ & $6$ \\
$2$ & $2$ & $3$ & $0$ & $1$ & $6$ & $7$ & $4$ & $5$ \\
$3$ & $3$ & $2$ & $1$ & $0$ & $7$ & $6$ & $5$ & $4$ \\
$4$ & $4$ & $5$ & $6$ & $7$ & $0$ & $1$ & $2$ & $3$ \\
$5$ & $5$ & $4$ & $7$ & $6$ & $1$ & $0$ & $3$ & $2$ \\
$6$ & $6$ & $7$ & $4$ & $5$ & $3$ & $2$ & $1$ & $0$ \\
$7$ & $7$ & $6$ & $5$ & $4$ & $2$ & $3$ & $0$ & $1$ 
\end{tabular}
\caption{Gyro table of $K(1)$}
\end{center}
\end{table}

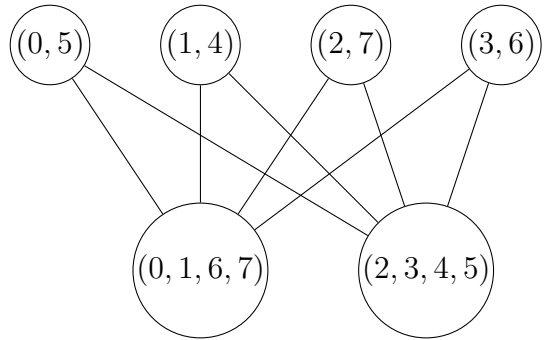
\begin{figure}[H]
\begin{center}
\begin{tikzpicture}
\node [circle,inner sep=.5pt,draw] (v1) at (0,0) {$(0,5)$}; \node [circle,inner sep=.5pt,draw] (v2) at (2,0) {$(1,4)$}; \node [circle,inner sep=.5pt,draw] (v3) at (4,0) {$(2,7)$}; \node [circle,inner sep=.5pt,draw] (v4) at (6,0) {$(3,6)$}; 
\node [circle,inner sep=.5pt,draw] (v5) at (2, -3) {$(0,1,6,7)$}; \node [circle,inner sep=.5pt,draw] (v6) at (5,-3) {$(2,3,4,5)$}; 
\draw (v1) -- (v5); \draw (v1) -- (v6); \draw (v2) -- (v5); \draw (v2) -- (v6); \draw (v3) -- (v5); \draw (v3) -- (v6); \draw (v4) -- (v5); 
\draw (v4) -- (v6);
\end{tikzpicture}
\caption{$\Phi(K(1), S= \{ 5,7\})$.}
\end{center}
\end{figure}
 
 \begin{table}[H]
 
 \begin{center} 
\begin{tabular}{c|cccccccc}
${\oplus}_{L(1)}$ & $0$ & $1$ & $2$ & $3$ & $4$ & $5$ & $6$ & $7$\\
\hline
$0$ & $0$ & $1$ & $2$ & $3$ & $4$ & $5$ & $6$ & $7$ \\
$1$ & $1$ & $0$ & $3$ & $2$ & $5$ & $4$ & $7$ & $6$ \\
$2$ & $2$ & $3$ & $0$ & $1$ & $6$ & $7$ & $4$ & $5$ \\
$3$ & $3$ & $2$ & $1$ & $0$ & $7$ & $6$ & $5$ & $4$ \\
$4$ & $4$ & $5$ & $6$ & $7$ & $0$ & $1$ & $2$ & $3$ \\
$5$ & $5$ & $4$ & $7$ & $6$ & $1$ & $0$ & $3$ & $2$ \\
$6$ & $6$ & $7$ & $5$ & $4$ & $3$ & $2$ & $0$ & $1$ \\
$7$ & $7$ & $6$ & $4$ & $5$ & $2$ & $3$ & $1$ & $0$ \\
\end{tabular}
\caption{Cayley table of $L(1)$}
\end{center}
\end{table}

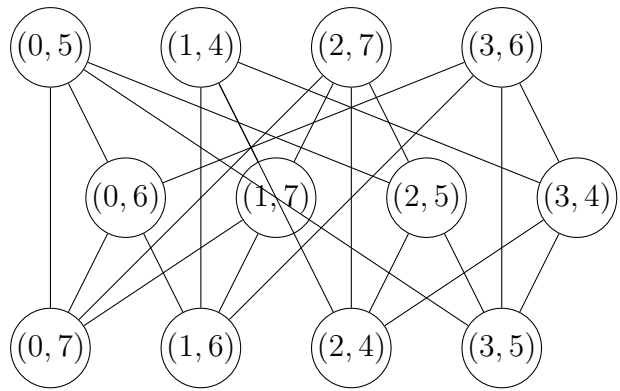
\begin{figure}[H]
\begin{center} 
\begin{tikzpicture}
\node [circle,inner sep=.5pt,draw] (v1) at (0,0) {$(0,5)$}; 
\node [circle,inner sep=.5pt,draw] (v2) at (2,0) {$(1,4)$}; 
\node [circle,inner sep=.5pt,draw] (v3) at (4,0) {$(2,7)$}; 
\node [circle,inner sep=.5pt,draw] (v4) at (6,0) {$(3,6)$}; 
\node [circle,inner sep=.5pt,draw] (v5) at (1,-2) {$(0,6)$}; 
\node[circle,inner sep=.5pt,draw]  (v6) at (3,-2) {$(1,7)$}; 
\node[circle,inner sep=.5pt,draw]  (v7) at (5, -2) {$(2,5)$}; 
\node[circle,inner sep=.5pt,draw]  (v8) at ( 7,-2) {$(3,4)$}; 
\node [circle,inner sep=.5pt,draw] (v9) at (0,-4) {$(0,7)$}; 
\node [circle,inner sep=.5pt,draw] (v10) at (2,-4) {$(1,6)$}; 
\node [circle,inner sep=.5pt,draw] (v11) at (4, -4) {$(2,4)$}; 
\node [circle,inner sep=.5pt,draw] (v12) at ( 6,-4) {$(3,5)$}; 
 \draw (v1) -- (v5); \draw (v1) -- (v7); \draw (v1) -- (v9); \draw (v1) -- (v12); 
\draw (v2) -- (v6); \draw (v2) -- (v8); \draw (v2) -- (v10); \draw (v2) -- (v11); 
\draw (v3) -- (v6); \draw (v3) -- (v7); \draw (v3) -- (v9); \draw (v3) -- (v11); 
\draw (v4) -- (v5); \draw (v4) -- (v8); \draw (v4) -- (v10); \draw (v4) -- (v12); 
\draw (v5) -- (v9); \draw (v5) -- (v10); \draw (v6) -- (v9); \draw (v6) -- (v10); 
\draw (v7) -- (v11); \draw (v7) -- (v12); \draw (v8) -- (v11); \draw (v8) -- (v12); 
\end{tikzpicture}
\caption{$\Phi(L(1), S= \{ 5,6,7\})$.}
\end{center}
\end{figure}

 \begin{table}[H]\begin{center}
\begin{tabular}{c|cccccccc}
${\oplus}_{M(1)}$ & $0$ & $1$ & $2$ & $3$ & $4$ & $5$ & $6$ & $7$\\
\hline
$0$ & $0$ & $1$ & $2$ & $3$ & $4$ & $5$ & $6$ & $7$ \\
$1$ & $1$ & $0$ & $3$ & $2$ & $5$ & $4$ & $7$ & $6$ \\
$2$ & $2$ & $3$ & $0$ & $1$ & $6$ & $7$ & $4$ & $5$ \\
$3$ & $3$ & $2$ & $1$ & $0$ & $7$ & $6$ & $5$ & $4$ \\
$4$ & $4$ & $5$ & $6$ & $7$ & $1$ & $0$ & $3$ & $2$ \\
$5$ & $5$ & $4$ & $7$ & $6$ & $0$ & $1$ & $2$ & $3$ \\
$6$ & $6$ & $7$ & $5$ & $4$ & $2$ & $3$ & $1$ & $0$ \\
$7$ & $7$ & $6$ & $4$ & $5$ & $3$ & $2$ & $0$ & $1$ \\
\end{tabular}
\caption{Gyro table of $M(1)$}
\end{center}
\end{table}

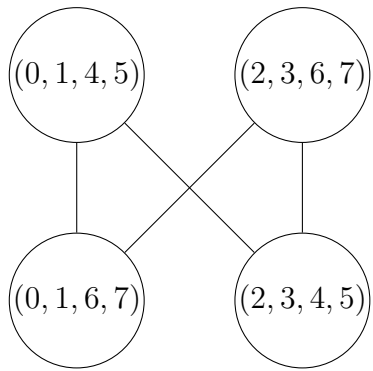
\begin{figure}[H]
  \begin{center}
\begin{tikzpicture}
\node [circle,inner sep=.5pt,draw] (v1) at (0,0) {$(0, 1, 4, 5)$}; \node [circle,inner sep=.5pt,draw] (v2) at (3,0) {$(2,3,6,7)$}; 
\node [circle,inner sep=.5pt,draw] (v3) at (0,-3) {$(0,1,6,7)$}; \node [circle,inner sep=.5pt,draw] (v4) at (3,-3) {$(2,3,4,5)$}; 
\draw (v1) -- (v3); \draw (v1) -- (v4); \draw (v2) -- (v3); \draw (v2) -- (v4); 
\end{tikzpicture}
\caption{$\Phi(M(1), S= \{ 5,7\})$.}
\end{center}
\end{figure}

\begin{table}\begin{center}
\begin{tabular}{c|cccccccc}
${\oplus}_{N(1)}$ & $0$ & $1$ & $2$ & $3$ & $4$ & $5$ & $6$ & $7$\\
\hline
$0$ & $0$ & $1$ & $2$ & $3$ & $4$ & $5$ & $6$ & $7$ \\
$1$ & $1$ & $0$ & $3$ & $2$ & $5$ & $4$ & $7$ & $6$ \\
$2$ & $2$ & $3$ & $0$ & $1$ & $6$ & $7$ & $4$ & $5$ \\
$3$ & $3$ & $2$ & $1$ & $0$ & $7$ & $6$ & $5$ & $4$ \\
$4$ & $4$ & $5$ & $6$ & $7$ & $1$ & $0$ & $3$ & $2$ \\
$5$ & $5$ & $4$ & $7$ & $6$ & $0$ & $1$ & $2$ & $3$ \\
$6$ & $6$ & $7$ & $5$ & $4$ & $2$ & $3$ & $1$ & $0$ \\
$7$ & $7$ & $6$ & $4$ & $5$ & $3$ & $2$ & $0$ & $1$ \\
\end{tabular}
\caption{Gyro table of $N(1)$}\end{center}
\end{table}
 
 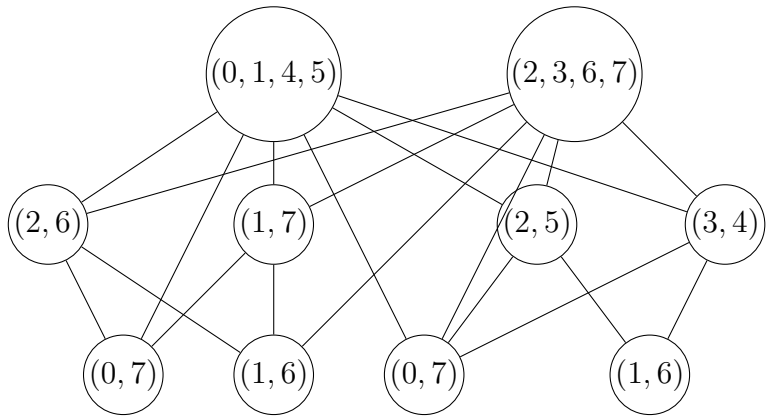
\begin{figure}[H] 
\begin{center} 
\begin{tikzpicture}
\node [circle,inner sep=.5pt,draw] (v1) at (0,0) {$(0,1,4,5)$}; 
\node [circle,inner sep=.5pt,draw] (v2) at (4,0) {$(2,3,6,7)$}; 
\node[circle,inner sep=.5pt,draw]  (v3) at (-3,-2) {$(2,6)$}; 
\node [circle,inner sep=.5pt,draw] (v4) at (0,-2) {$(1,7)$}; 
\node [circle,inner sep=.5pt,draw] (v5) at (3.5,-2) {$(2,5)$}; 
\node [circle,inner sep=.5pt,draw] (v6) at (6,-2) {$(3,4)$};
 \node [circle,inner sep=.5pt,draw] (v7) at (-2, -4) {$(0,7)$}; 
 \node [circle,inner sep=.5pt,draw] (v8) at ( 0,-4) {$(1,6)$}; 
 \node [circle,inner sep=.5pt,draw] (v9) at (2,-4) {$(0,7)$}; 
 \node [circle,inner sep=.5pt,draw] (v10) at (5,-4) {$(1,6)$}; 
 \draw (v1) -- (v3); \draw (v1) -- (v4); \draw (v1) -- (v5); \draw (v1) -- (v6); \draw (v1) -- (v7); \draw(v1) -- (v9); 
\draw (v2) -- (v3); \draw (v2) -- (v4); \draw (v2) -- (v5); \draw (v2) -- (v6);  \draw (v2) -- (v9); \draw (v2) -- (v8);
\draw  (v3) -- (v7); \draw (v3) -- (v8); 
\draw  (v4) -- (v7); \draw (v4) -- (v8); 
\draw (v5) -- (v9); \draw (v5) -- (v10); 
\draw (v6) -- (v9); \draw (v6) -- (v10); 
\end{tikzpicture}
\caption{ $\Phi(N(1), S= \{ 5,6,7\})$.}
\end{center}
\end{figure}

 \begin{table} \begin{center}
\begin{tabular}{c|cccccccc}
${\oplus}_{O(1)}$ & $0$ & $1$ & $2$ & $3$ & $4$ & $5$ & $6$ & $7$\\
\hline
$0$ & $0$ & $1$ & $2$ & $3$ & $4$ & $5$ & $6$ & $7$ \\
$1$ & $1$ & $0$ & $3$ & $2$ & $5$ & $4$ & $7$ & $6$ \\
$2$ & $2$ & $3$ & $0$ & $1$ & $6$ & $7$ & $4$ & $5$ \\
$3$ & $3$ & $2$ & $1$ & $0$ & $7$ & $6$ & $5$ & $4$ \\
$4$ & $4$ & $5$ & $7$ & $6$ & $1$ & $0$ & $2$ & $3$ \\
$5$ & $5$ & $4$ & $6$ & $7$ & $0$ & $1$ & $3$ & $2$ \\
$6$ & $6$ & $7$ & $5$ & $4$ & $2$ & $3$ & $1$ & $0$ \\
$7$ & $7$ & $6$ & $4$ & $5$ & $3$ & $2$ & $0$ & $1$ \\
\end{tabular}
\caption{Gyro table of $O(1)$}\end{center}  
\end{table}

\tikzset{me/.style={to path={
\pgfextra{%
 \pgfmathsetmacro{\startf}{-(#1-1)/2}  
 \pgfmathsetmacro{\endf}{-\startf} 
 \pgfmathsetmacro{\stepf}{\startf+1}}
 \ifnum 1=#1 -- (\tikztotarget)  \else
     let \p{mid}=($(\tikztostart)!0.5!(\tikztotarget)$) 
         in
\foreach \i in {\startf,\stepf,...,\endf}
    {%
     (\tikztostart) .. controls ($ (\p{mid})!\i*6pt!90:(\tikztotarget) $) .. (\tikztotarget)
      }
      \fi   
     \tikztonodes
}}}     

\begin{figure}[H]  
 \begin{center}
\begin{tikzpicture}
\node [circle,inner sep=.5pt,draw] (v1) at (0,0) {$(0, 1, 4, 5)$}; \node [circle,inner sep=.5pt,draw] (v2) at (3,0) {$(2,3,6,7)$}; 
\node [circle,inner sep=.5pt,draw] (v3) at (0,-3) {$(0,1,6,7)$}; \node [circle,inner sep=.5pt,draw] (v4) at (3,-3) {$(2,3,5,4)$}; 
\draw (v1) edge[me=2] (v3); 
\draw (v1) edge[me=2] (v4); \draw (v2) edge[me=2] (v3); \draw (v2) edge[me=2] (v4); 
\end{tikzpicture}
\caption{$\Phi(O(1), S= \{ 5,7\})$.}
\end{center}
\end{figure}
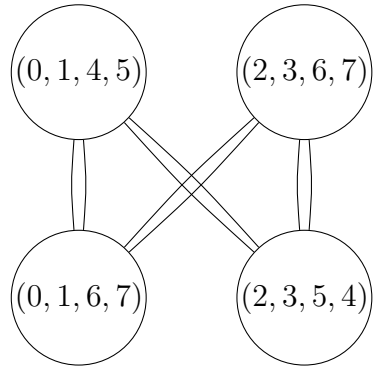

Note that the set of generators mentioned for the above gyro-groups is the smallest possible generator set.

\begin{exa}
Consider two gyro-groups $ K(1) $ and $ G_{8,2} $ with their gyro tables. These two gyro-groups are isomorphic to each other. Now we compute their Cayley graphs.
\begin{table}
\begin{center}
\footnotesize{
\begin{tabular}{c|llllllll lc|llllllll}
$\oplus_{K(1)}$ &  $0$ & $1$ & $2$ & $3$ & $4$ & $5$ & $6$ & $7$ & & $\oplus$ &  $0$ & $1$ & $2$ & $3$ & $4$ & $5$ & $6$ & $7$ \\
\hline 
$0$ & $0$ & $1$ & $2$ & $3$ & $4$ & $5$ & $6$ & $7$ & & $0$ & $0$ & $1$ & $2$ & $3$ & $4$ & $5$ & $6$ & $7$\\
$1$ & $1$ & $0$ & $3$ & $2$ & $5$ & $4$ & $7$ & $6$ &  & $1$ & $1$ & $7$ & $6$ & $0$ & $5$ & $2$ & $4$ & $3$\\
$2$ &  $2$ & $3$ & $ 0$ & $1$ & $6$ & $7$ & $4$ & $5$ & & $2$ &  $2$ & $5$ & $7$ & $6$ & $0$ & $3$ & $1$ & $4$ \\
$3$ &  $3$ & $2$ & $1$ & $0$ & $7$ &  $6$ & $5$ & $4$  & & $3$ &  $3$ & $0$ & $5$ & $7$ & $6$ & $4$ & $2$ & $1$ \\
$4$ &  $4$ & $5$ & $6$ & $7$ & $0$ & $1$ & $2$ & $3$  &  & $4$ &  $4$ & $6$ & $0$ & $5$ & $7$ & $1$ & $3$ & $2$ \\
$5$ &  $5$ & $4$ & $7$ & $6$ & $1$ & $0$ & $3$ & $2$  &  & $5$ &  $5$ & $2$ & $3$ & $4$ & $1$ & $7$ & $0$ & $6$ \\
$6$ &  $6$ & $7$ & $4$ & $5$ & $3$ & $2$ & $1$ & $0$  &  & $6$ &  $6$ & $4$ & $1$ & $2$ & $3$ & $0$ & $7$ & $5$ \\
$7$ &  $7$ & $6$ & $5$ & $4$ & $2$ & $3$ & $0$ & $1$   & & $7$ &  $7$ & $3$ & $4$ & $1$ & $2$ & $6$ & $5$ & $0$ 
\end{tabular}
}
\caption{ Cayley table $K(1)$ and the gyro table of $G_8$}\end{center}
\end{table}
Their corresponding Cayley graphs are respectively presented in two following figures, which you can see that 
 $$Cay(K(1), S=\{5,7\}) \ncong Cay(G_{8,2}, S=\{7,8\})$$

\tikzset{me/.style={to path={
\pgfextra{%
 \pgfmathsetmacro{\startf}{-(#1-1)/2}  
 \pgfmathsetmacro{\endf}{-\startf} 
 \pgfmathsetmacro{\stepf}{\startf+1}}
 \ifnum 1=#1 -- (\tikztotarget)  \else
     let \p{mid}=($(\tikztostart)!0.5!(\tikztotarget)$) 
         in
\foreach \i in {\startf,\stepf,...,\endf}
    {%
     (\tikztostart) .. controls ($ (\p{mid})!\i*6pt!90:(\tikztotarget) $) .. (\tikztotarget)
      }
      \fi   
     \tikztonodes
}}}   
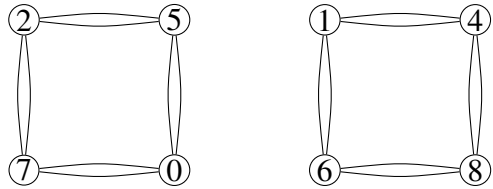
\begin{figure}[H]
\begin{center}
\begin{tikzpicture}
\node[circle,inner sep=.5pt,draw] (v1) at (0,0) {2};
\node[circle,inner sep=.5pt,draw] (v2) at (2,0) {5};
\node[circle,inner sep=.5pt,draw] (v3) at (0,-2) {7};
\node[circle,inner sep=.5pt,draw] (v4) at (2,-2) {0};
\node[circle,inner sep=.5pt,draw] (v5) at (4,0) {1};
\node[circle,inner sep=.5pt,draw] (v6) at (6,0) {4};
\node[circle,inner sep=.5pt,draw] (v7) at (4,-2) {6};
\node[circle,inner sep=.5pt,draw] (v8) at (6,-2) {8};
\draw (v1) edge[me=2] (v2); 
\draw (v1) edge[me=2] (v3); 
 \draw (v2) edge[me=2] (v4); 
 \draw (v3) edge[me=2] (v4); 
 \draw (v5) edge[me=2] (v6); 
\draw (v5) edge[me=2] (v7); 
 \draw (v6) edge[me=2] (v8); 
 \draw (v7) edge[me=2] (v8);  
\end{tikzpicture} 
\caption{$Cay(K(1), S= \{ 5,7\})$.}
 \end{center}
 \end{figure}

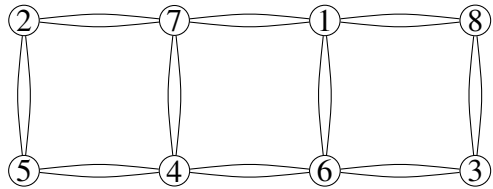
\begin{figure}[H] 
 \begin{center}
\begin{tikzpicture}
\node[circle,inner sep=.5pt,draw] (v1) at (0,0) {2};
\node[circle,inner sep=.5pt,draw] (v2) at (2,0) {7};
\node[circle,inner sep=.5pt,draw] (v3) at (0,-2) {5};
\node[circle,inner sep=.5pt,draw] (v4) at (2,-2) {4};
\node[circle,inner sep=.5pt,draw] (v5) at (4,0) {1};
\node[circle,inner sep=.5pt,draw] (v6) at (6,0) {8};
\node[circle,inner sep=.5pt,draw] (v7) at (4,-2) {6};
\node[circle,inner sep=.5pt,draw] (v8) at (6,-2) {3};
\draw (v1) edge[me=2] (v2); 
\draw (v1) edge[me=2] (v3); 
 \draw (v2) edge[me=2] (v4); 
 \draw (v3) edge[me=2] (v4); 
 \draw (v2) edge[me=2] (v5); 
\draw (v4) edge[me=2] (v7);  
\draw (v5) edge[me=2] (v6); 
\draw (v5) edge[me=2] (v7); 
 \draw (v6) edge[me=2] (v8); 
 \draw (v7) edge[me=2] (v8);  
\end{tikzpicture} 
\caption{ $Cay(G_{8,2}, S= \{ 7, 8\})$.}
 \end{center}
 \end{figure}
  
\end{exa}

Here there is an example of two isomorphic gyro-groups with non-isomorphic Cayley graphs.
\begin{exa}
Consider two isomorphic gyro-groups $ N(1) $ and $ G_{8,6} $. We  obtain the structure  of $\Phi( G_{8,6} ,  S= \lbrace 7, 8\rbrace) $.

\begin{figure}[H]
\begin{center} 
\begin{tikzpicture}
\node [circle,inner sep=.5pt,draw] (v1) at (0,0) {$(1,3,5,7)$}; \node [circle,inner sep=.5pt,draw] (v2) at (4,0) {$(2,4,6,8)$}; 
\node [circle,inner sep=.5pt,draw] (v3) at (-3,-2) {$(1,8)$}; \node [circle,inner sep=.5pt,draw] (v4) at (0,-2) {$(2,5)$}; 
\node[circle,inner sep=.5pt,draw]  (v5) at (3.5, -2) {$(3,6)$}; \node [circle,inner sep=.5pt,draw] (v6) at (6,-2) {$(4,7)$};
 \draw (v1) -- (v3); \draw (v1) -- (v4); \draw (v1) -- (v5); \draw (v1) -- (v6); 
 \draw (v2) -- (v3); \draw (v2) -- (v4); \draw (v2) -- (v5); \draw (v2) -- (v6);  
 \end{tikzpicture}
\caption{$\Phi(G_{8,6}, S=\{ 7,8\})$}
\end{center}
\end{figure}
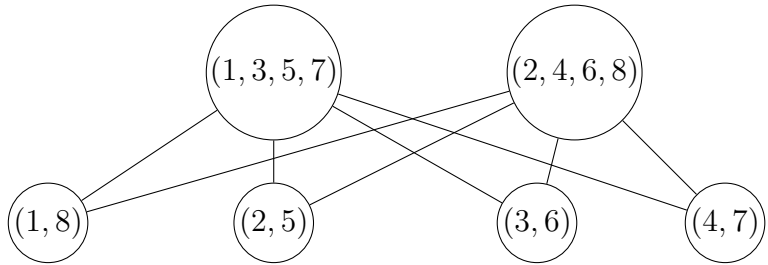

You can see that
$$\Phi( G_{8,6}, \{7,8\}) \ncong \Phi(N(1), \{5,6,7\})$$ 
\end{exa}

There is an example that shows there are isomorphic gyro-groups with non-isomorphic $G$-graphs.
In the following, we will examine the $G$-graphs of a special gyro-group, which is called the gyro-commutative gyro-group.

\begin{defn} \cite{Suksumran2016}(Gyro-commutative low)
A gyro-group $ G $   is called a gyro-commutative gyro-group, if for all $ a, b \in G  $,   we have:
$$a \oplus b =gyr[a, b]( b \oplus a ) $$
\end{defn}

\begin{exa} 
 Consider two gyro-groups Isomorphic $ M(1) $ and $ G_{8,5} $. We know $ G_{8,5}  $ is gyro-commutative.
We compute $G$-graph of $ G_{8,5} $. We obtain the following graph:

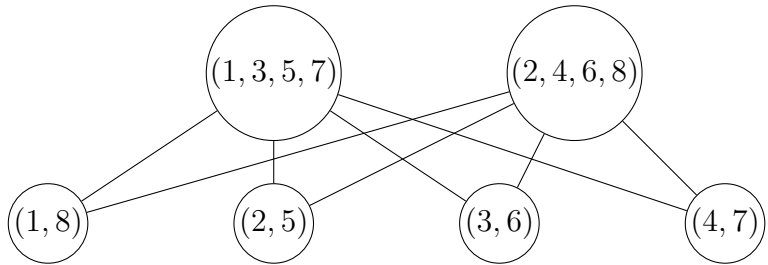
\begin{figure}[H]
\begin{center}
\begin{tikzpicture}
\node [circle,inner sep=.5pt,draw] (v1) at (0,0) {$(1,3,5,7)$}; \node [circle,inner sep=.5pt,draw] (v2) at (4,0) {$(2,4,6,8)$}; 
\node [circle,inner sep=.5pt,draw] (v3) at (-3,-2) {$(1,8)$}; \node [circle,inner sep=.5pt,draw] (v4) at (0,-2) {$(2,5)$};  
\node [circle,inner sep=.5pt,draw] (v5) at (3,-2) {$(3,6)$};  \node[circle,inner sep=.5pt,draw]  (v6) at (6,-2) {$(4,7)$};  
\draw (v1) -- (v3); \draw (v1) -- (v4); \draw (v1) -- (v5); \draw (v1) -- (v6); 
\draw (v2) -- (v3); \draw (v2) -- (v4); \draw (v2) -- (v5); \draw (v2) -- (v6); 
\end{tikzpicture}
\caption{ $\Phi(G_{8,5}, S= \{ 7,8\})$.}
\end{center}
\end{figure}

You see that the $G$-graphs of these two isomorphic  gyro-commutative gyro-groups are not isomorphic to each other.
\end{exa}

\begin{rem}
The structure of the $G$-graph of the gyro-groups $K(1)$, $ M(1)$ and $ O(1)$ are as follows:
\begin{itemize}
\item[I.] $\Phi( K(1),S= \lbrace 5, 7\rbrace  ) $  is a complete bipartite graph isomorphic to $ K_{2, 4} $
\item[II. ]  $ \Phi ( M(1), S= \lbrace 5, 7\rbrace  ) $  is a complete bipartite graph isomorphic to $ K_{2, 2} $ With a double edge.
\item[III.] $\Phi ( O(1), S= \lbrace 5, 7\rbrace  ) $  is complete bipartite graph $ K_{2, 2} $ With a double edge.
\end{itemize}
\end{rem}

 \subsection{Cayley graph and $G$-graph of $2$-gyro-groups}
 
 Ashrafi, et. al, In \cite{Ashrafi} have presented an algebraic approach to constructing new gyro-groups. They  
 constructed a gyro-group $ G(n) $, $ n \geq 3 $, of order $ 2 ^{n} $ by considering a cyclic group of order
$ 2^{n}$ . It is proved that all proper sub-gyro-groups of $ G(n) $ are either cyclic or dihedral groups.
 We will examine the $G$-graph of $ 2 $-gyro-groups. 
 
\begin{defn} \cite{Ashrafi}
A class of $2$-gyro-groups constructed by Ashrafi, et. al, in such a way that:
\begin{center}
$ P(n) = \lbrace 0, 1, 2,\cdots , 2^{n -1} - 1 \rbrace$, $ H(n) =  \lbrace 2^{n - 1}, 2^{n-1} + 1, \cdots , 2^{n} - 1\rbrace $,
\end{center}
 where $ G(n) = P(n) \cup H(n)$, for a natural number  $ n \geq 3 $,  that  $ P(n) $ is a cyclic group under addition modulo $ m = 2^{n - 1}$ and  $ H(n) = P(n) + m $. This shows that
$ G(n) = P(n) \cup (P(n) + m)$.
\end{defn}
You can refer to \cite{Ashrafi} for more details. 

\begin{thm}\cite{Ashrafi}
$(G(n), \oplus)$ is a gyro-group.
\end{thm}

\begin{exa}\cite{Ashrafi}
Consider the gyro-group 
$$G(4) = ( \lbrace 0, 2, 3, 4, 5, 6, 7, 8, 9, 10, 11, 12, 13, 14, 15 \rbrace, \oplus ),$$ 
with $A = (8, 12)(9, 13)(10, 14)(11, 15)$.
The Cayley table and the $G$-graph  $ \Phi (G(4), H(4)) $ are presented in the following.
\begin{table}[H]\begin{center}
{\footnotesize{
\begin{tabular}{c|ccccccccccccccccc}
        $\oplus$ & 0 & 1 & 2 & 3 & 4 & 5 & 6 & 7 & 8 & 9 & 10 & 11 & 12 & 13 & 14 & 15  \\ 
        \hline
		0 & 0 & 1 & 2 & 3 & 4 & 5 & 6 & 7 & 8 & 9 & 10 & 11 & 12 & 13 & 14 & 15  \\
		1 & 1 & 2 & 3 & 4 & 5 & 0 & 7 & 9 & 10 & 11 & 8 & 12 & 13 & 14 & 15 & 8  \\
		2 & 2 & 3 & 4 & 5& 6 & 7 & 0& 1&  10 & 11 & 12 & 13 & 14 & 15 & 8 & 9  \\
		3 & 3 & 4 & 5 & 6 & 7& 0 & 1 & 2 & 11 & 12 & 13 & 14 & 15 & 8 & 9 & 10  \\
		4 & 4 & 5 & 6 & 7 & 0 & 1 & 2 & 3 & 12 & 13 & 14 & 15 & 8 & 9 & 10 & 11 \\
		5 & 5 & 6 & 7 & 0& 1 & 2 & 3 & 4 & 13 & 14 & 15 & 8 & 9 & 10 & 11 & 12  \\
		6 & 6 & 7 & 0 & 1 & 2 & 3 & 4 & 5 & 14 & 15 & 8 & 9 & 10 & 11 & 12 & 13  \\
		7 & 7 & 0 & 1 & 2 & 3 & 4 & 5 & 6 & 15 & 8 & 9 & 10 & 11 & 12 & 13 & 14   \\
		8 & 8 & 11& 14 & 9 & 12 & 15 & 10 & 13 & 0 & 3 & 6 & 1 & 4 & 7 & 2 & 5  \\
		9 & 9 & 12 & 15 & 10 & 13 & 8 & 11 & 14 &  5 & 0 & 3 & 6 & 1 & 4 & 7 & 2     \\
		10 & 10 & 13 & 8 & 11 & 14 & 9 & 12 & 15 & 2 & 5 & 0 & 3 & 6 & 1& 4 & 7  \\
		11 & 11 & 14 & 9 & 12 & 15 & 10 & 13 & 8 & 7 & 2 & 5 & 0 & 3 & 6 & 1 & 4 \\
		12 & 12 & 15 & 10 & 13 & 8 & 11 & 14 & 9 & 4 & 7 & 2 & 5 & 0 & 3 & 6 & 1  \\
		13 & 13 & 8 & 11 & 14 & 9 & 12 & 15 & 10 & 1 & 4 & 7 & 2 & 5 & 0 & 3 & 6 \\
		14 & 14 & 9 & 12 & 15 & 10 & 13 & 8 & 11 & 6 & 1 & 4 & 7 & 2 & 5 & 0 & 3  \\
		15 & 15 & 10 & 13 & 8 & 11 & 14 & 9 & 12 & 3 & 6 & 1 & 4 & 7 & 2 & 5 & 0 \\
\end{tabular}
}}
\caption{The Cayley Table of $G(4)$}\label{ta1}
\end{center}
\end{table}

We know that $ m=2^{n-1}= 2^{4-1}= 8 $ then
$$H(3)=\lbrace 0,1,2,\cdots, 2^{n-1} -1 \rbrace= \lbrace 0, 1, 2, 3, 4, 5, 6, 7\rbrace$$
and $  \Phi (G(4), H(4)) $ has $ 8 $ parts of   the vertex set:
 $$ V= V_{8} \bigcup  V_{9} \bigcup V_{10} \bigcup V_{11} \bigcup V_{12} \bigcup V_{13} \bigcup V_{14} \bigcup V_{15}.$$
 thus we have:
$$ V= \lbrace  ( 0, 1), (1, 11), (2, 14), (3, 9), (4, 12), (5, 15), (6, 10), (7, 13) \rbrace$$
$$\bigcup  \lbrace  ( 0, 9), (1, 12), (2, 15), (3, 10), (4, 13), (5, 8), (6, 11), (7, 1) \rbrace$$
$$\bigcup \lbrace ( 0, 10), (1, 13), (2, 8), (3, 11), (4, 14), (5, 9), (6, 12), (7, 15) \rbrace$$
$$\bigcup \lbrace  ( 0, 11), (1, 14), (2, 9), (3, 12), (4, 15), (5, 10), (6, 13), (7, 8) \rbrace$$
$$\bigcup  \lbrace ( 0, 12), (1, 15), (2, 10), (3, 13), (4, 8), (5, 11), (6, 14), (7, 9) \rbrace$$
$$\bigcup \lbrace  ( 0, 13), (1, 8), (2, 11), (3, 14), (4, 9), (5, 12), (6, 15), (7, 10) \rbrace$$
$$\bigcup \lbrace ( 0, 14), (1, 9), (2, 12), (3, 15), (4, 10), (5, 13), (6, 8), (7, 11) \rbrace$$ 
 $$\bigcup  \lbrace ( 0, 15), (1, 10), (2, 13), (3, 8), (4, 11), (5, 14), (6, 9), (7, 12)  \rbrace .$$

This graph is $2(m - 1) =2( 8 - 1) = 14 $-regular. Moreover for every    $ a \in G(4) $ and $ j \in  H(4)$, $ \langle j \rangle \cong \mathbb{Z}_{2} $ we can see that $gyr[a; j]( \langle j \rangle) = \langle j \rangle $, then $  \langle j \rangle $ is the $ L$-sub-gyro-group of $G(4)$, In this case,
 according to Lagrange's theorem for gyro-groups:
$$ \vert V_{8} \vert =  \vert V_{9} \vert = \vert V_{10} \vert = \vert V_{11} \vert = \vert V_{12} \vert = \vert V_{13} \vert = \vert V_{14} \vert = \vert V_{15} \vert = \vert G(4) : V_{j} \vert = 2 .$$
\end{exa}

One of the important concept in graph theory is hamiltonicity of a graph.  Here there is a theorem which says that the structure of the $G$-graph of the gyro-groups $(G(n), H(n))$ is hamiltonian.

\begin{thm}\cite{Farzaneh2022}
The $G$-gyro-graph $\Phi (G(n), H(n))$ is connected and hamiltonian.
\end{thm}
\begin{proof}
For every $ j \in  H(n) =  \lbrace m, m + 1, \cdots, 2m - 1 \rbrace$, $ m = 2^{n - 1}$, $ n \geq 3$ and $ a \in G(n) $, we have
\begin{center}
$gyr[a, s](\langle j \rangle) = \langle j \rangle $, that $ s \in \langle j \rangle $.
\end{center}
Consequently $ \langle j \rangle \cong \mathbb{Z}_{2} $ is an $ L $-sub-gyro-group of $ G(n) $. According to the establishment of Lagrange's theorem holds for $ L $-sub-gyro-groups, we have
$$ \vert G(n), \langle j \rangle \vert= \frac{2n}{2} = m,$$
and $ \Phi (G(n), H(n) )$ is an $m$- partite, each part with $ m $ vertices of
 $$ (j) \oplus x = (x, j \oplus x),~~~~ x \in G(n),$$
  that has at least two intersections with the vertices in other levels. Indeed the graph is $2(m - 1)$-regular.  Also each $V_{j} $ has a copy of $G(n) $ which induces that $\Phi (G(n), H(n))$, is connected. On the other hand,  there is a hamiltonian path in this graph. As a result, the graph is connected.
 \end{proof}

\begin{exa}
 We want to get $ \Phi (G(4), P(4))$.  Consider the subgroup
 
 $\left\langle  1 \right\rangle \cong  \mathbb{Z}_{m} = \mathbb{Z}_{8}$, for every $a \in G(4) $ and $ x \in \left\langle 1 \right\rangle $, we see that 
 
 $gyr[a; x](\left\langle 1\right\rangle =  \left\langle 1 \right\rangle $, then $\left\langle 1 \right\rangle  $ is an $ L $-subgyro-group of $ G(4) $ and

 $  \vert V_{1} \vert = \vert G(4) : \langle 1 \rangle \vert = 2$.
 \begin{flushleft}
$ V_{1} = \lbrace (1) \oplus x = \lbrace (0, 1, 2, 3, 4, 5, 6, 7), (8, 9, 10, 11, 12, 13, 14, 15) \rbrace $,

$V_{2} = \lbrace (0, 2, 4, 6), (1, 3, 5, 7), (8, 10, 12, 14), (9, 11, 13, 15) \rbrace$,

$V_{3} = \lbrace ( 0, 1, 2, 3, 4, 5, 6, 7),(8, 9, 10, 11, 12, 13, 14, 15) \rbrace $,

$V_{4} = \lbrace (0, 4), (1,5), (6, 2), (7, 3), (8, 12), (13, 9), (14, 10), (11, 15) \rbrace$,

$V_{5} = \lbrace (0, 1, 2, 3, 4, 5, 6, 7), (8, 9, 10, 11, 12, 13, 14, 15) \rbrace$,

$V_{6} = \lbrace  (0, 2, 4, 6), (1, 3, 5, 7), (8, 10, 12, 14), (9, 11, 13, 15) \rbrace $,

$V_{7} =  \lbrace (0, 1, 2, 3, 4, 5, 6, 7), (8, 9, 10, 11, 12, 13, 14, 15) \rbrace$.
\end{flushleft}

Since $  \Phi (G(n), H(n) ) $ hasn't any loops in $ G $, then based on   \cite{Badaoui}

$ d(v) = o(s)(k - 1)$ for all $ v \in V_{s} $, thus we have that $ \Phi (G(4), P(4))$ is a connected, $4$-partite  $G$-graph with
$deg(v_{1}) =56$, $deg(v_{2}) = 28 $, $deg(v_{3}) = 56 $, $deg(v_{4}) = 16 $, $ deg(v_{5}) = 56 $, $ deg(v_{6}) = 28 $, $ deg(v_{7}) = 56 $. 
\end{exa}

\subsection{$G$-graph of the gyro-commutative  gyro-groups}
Now we turn to the $G$-graph of the gyro-commutative  gyro-groups.

\begin{thm}\cite{Suksumran2015}
Let $p $ and  $ q$  be primes and let $ G $ be a gyro-group of order $ pq $. If $p= q$, then $G $ is a group. If $p \neq q$, then $ G $ is generated by two elements; one has order $p $ and the other has order $ q $.
\end{thm}
\begin{lem}\cite{Ashrafi2022}
There is a unique gyro-group of order $15 $ which is gyro-commutative
\end{lem}

 In the following we bring the structure of the $G$-graph of $G_{15}$. 
\begin{exa}
Consider the gyro-group $ G_{15} $. According to previous theorem and lemma, $ G_{15} $ is a gyro-commutative gyro-group and $ G_{15}= \langle 1, 4 \rangle $. 

\begin{table}[H]
	\centering
	{\footnotesize{
			\begin{tabular}{c|ccccccccccccccc}
				$\oplus$ & 0 & 1 & 2 & 3 & 4 & 5 & 6 & 7 & 8 & 9 & 10 & 11 & 12 & 13 & 14\\ 
				\hline
				0  & 0  & 1  & 2  & 3  & 4  & 5  & 6  & 7  & 8  & 9  & 10 & 11 & 12 & 13 & 14 \\ 
				1  & 1  & 2  & 0  & 4  & 6  & 11 & 3  & 14 & 13 & 7  & 8  & 12 & 5  & 10 & 9 \\ 
				2  & 2  & 0  & 1  & 6  & 3  & 12 & 4  & 9  & 10 & 14 & 13 & 5  & 11 & 8  & 7 \\ 
				3  & 3  & 4  & 5  & 7  & 8  & 9  & 13 & 0  & 1  & 2  & 12 & 6  & 14 & 11 & 10\\ 
				4  & 4  & 10 & 8  & 11 & 13 & 1  & 5  & 6  & 14 & 0  & 7  & 2  & 9  & 12 & 3\\ 
				5  & 5  & 14 & 12 & 9  & 7  & 8  & 2  & 11 & 0  & 10 & 3  & 4  & 6  & 1  & 13\\ 
				6  & 6  & 11 & 4  & 13 & 10 & 3  & 14 & 8  & 12 & 1  & 2  & 9  & 7  & 5  & 0\\ 
				7  & 7  & 8  & 9  & 0  & 1  & 2  & 11 & 3  & 4  & 5  & 14 & 13 & 10 & 6  & 12\\ 
				8  & 8  & 13 & 6  & 10 & 11 & 0  & 12 & 4  & 5  & 3  & 9  & 7  & 2  & 14 & 1\\ 
				9  & 9  & 5  & 11 & 14 & 0  & 6  & 7  & 10 & 2  & 12 & 1  & 3  & 13 & 4  & 8\\ 
				10 & 10 & 3  & 13 & 12 & 5  & 14 & 8  & 2  & 9  & 6  & 11 & 0  & 1  &  7 & 4\\ 
				11 & 11 & 12 & 7  & 1  & 14 & 4  & 9  & 13 & 6  & 8  & 0  & 10 & 3  &  2 & 5\\ 
				12 & 12 & 6  & 3  & 8  & 9  & 7  & 10 & 1  & 11 & 13 & 5  & 14 & 4  & 0  & 2\\ 
				13 & 13 & 7  & 14 & 2  & 12 & 10 & 1  & 5  & 3  & 4  & 6  & 8  & 0  & 9  & 11\\ 
				14 & 14 & 9  & 10 & 5  & 2  & 13 & 0  & 12 & 7  & 11 & 4  & 1  & 8  & 3  & 6 \\ 
				
			\end{tabular}
	}}
	\caption{Cayley table for the gyro-group $G_{15}$.} \label{tab: operation G15}
\end{table}
Since $ o(1)= 3 $ and $ o(4) = 5 $, then
$$V_{1}=\lbrace  (0,1,2), (3,4,6), (5,11,12), (7,9,14), (8, 10, 13) \rbrace $$
$$V_{4}= \lbrace (0,4,9,12,13), (1, 5, 6, 7, 10), (2,3,8,11,14)\rbrace, $$
and if we draw its associated $G$-graph, we  will have:

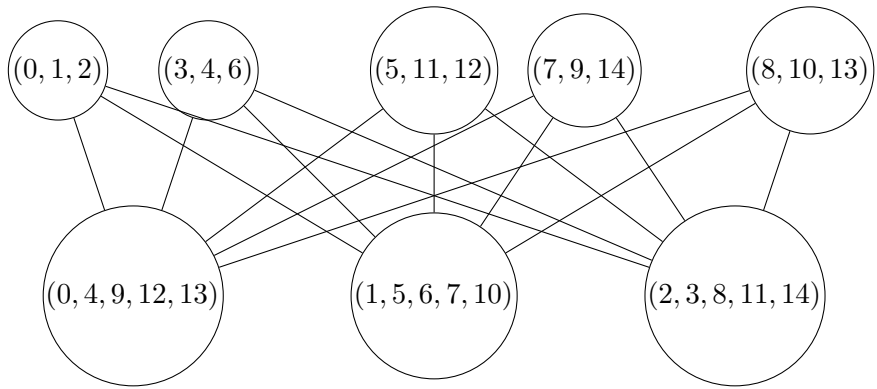
\begin{figure}[H]
\begin{center} 
\small{
\begin{tikzpicture}
\node [circle,inner sep=.5pt,draw] (v1) at (-7,0) {$(0,1,2)$}; 
\node [circle,inner sep=.5pt,draw] (v2) at (-5,0) {$(3,4,6)$};  \node [circle,inner sep=.5pt,draw] (v3) at (-2,0) {$(5,11,12)$}; 
\node [circle,inner sep=.5pt,draw] (v4) at (0,0) {$(7,9,14)$}; \node [circle,inner sep=.5pt,draw] (v5) at (3,0) {$(8,10,13)$}; 
\node [circle,inner sep=.5pt,draw] (v6) at (-6,-3) {$(0,4,9,12,13)$}; \node [circle,inner sep=.5pt,draw] (v7) at (-2, -3) {$(1,5,6,7,10)$}; 
\node [circle,inner sep=.5pt,draw] (v8) at ( 2,-3) {$(2,3,8,11,14)$}; 
 \draw (v1) -- (v6); \draw (v1) -- (v7); \draw (v1) -- (v8);
  \draw (v2) -- (v6); \draw (v2) -- (v7); \draw(v2) -- (v8); 
\draw (v3) -- (v6); \draw (v3) -- (v7); \draw (v3) -- (v8); 
\draw (v4) -- (v6);  \draw (v4) -- (v7); \draw (v4) -- (v8);
\draw (v5) -- (v6); \draw (v5) -- (v7); \draw (v5) -- (v8); 
\end{tikzpicture}
}
\caption{ $\Phi(G_{15}, S= \{ 1,4\})$.}
\end{center}
\end{figure}

\end{exa}

In gyro-groups including $ G(n) $, we have, if  $ S $ is a generating set, the graph is not connected, (In fact, it should be left generating) Now, the question is raised, is it possible that gyro-commutative gyro-group do not apply in this condition, and if $S $ is the generator set (and not the left generator), is the graph connected?
We are currently investigating other gyro-commutative gyro-group, one of which is the gyro-commutative gyro-group of order 21. According to reference \cite{Ashrafi2022}, we have only one gyro-group of this order, which is gyro-commutative gyro-group.

\subsection{$G$-graph of the  dihedralized gyro-groups}

In \cite{Maungchang2022} Dihedralized gyro-groups are introduced.  Maungchang and et.al, in \cite{Maungchang2022} have shown  with a gyro-commutative gyro-group satisfying the skew left loop property in hand, a new gyro-group can be constructed.

\begin{defn}\cite{Maungchang2022}
A gyro-group $G$ is said to have the skew left loop property if for all $ a, b \in G $
$$gyr [a \oplus b, \ominus b]=gyr[a, b] . $$ 
\end{defn}

\begin{defn}\cite{Maungchang2022}
 A gyro-group $ G $ is dihedralizable if it is gyro-commutative and has the skew left
loop property.
\end{defn}

In this article, we try to draw the $G$-graphs of these gyro-groups. These gyro-groups typically have a generating set of $ S=\lbrace (a, 0), (b, 0), (e, 1)\rbrace $.

\begin{exa}
Consider the gyro-group $ Dih(G_{8}) $ with its the gyration table.

\begin{table}
\begin{center}{\small{
\begin{tabular}{c|cccccccc}
$\oplus$& $(0,0)$ & $(1,0)$ & $(2,0)$ & $(3,0)$ & $(4,0)$ & $(5,0)$ & $(6,0)$ & $(7,0)$ \\ 
\hline \\
$ (0,0)$ & $(0,0)$ & $ (1,0)$ & $(2,0)$ & $(3,0)$ & $(4,0)$ & $(5,0)$	& $(6,0)$ & $(7,0)$  \\ 
$(1,0)$ & $(1,0)$ &  $(3,0)$ &	$(0,0)$ & $(2,0)$ & $(7,0)$ & $(4,0)$ & $(5,0)$ & $(6,0)$ \\
 $(2,0)$ & $(2,0)$ & $(0,0)$ &	$(3,0)$ & $(1,0)$ & $(5,0)$ & $(6,0)$ & $(7,0)$ & $(4,0)$ \\
(3,0)	&(3,0)&(2,0)&(1,0)	&(0,0)&(6,0)&(7,0)&(4,0)&(5,0) \\
(4,0)	&(4,0)&(5,0)&(7,0)&(6,0)&(3,0)&(2,0)&(0,0)&(1,0)\\
(5,0)	&(5,0)&(6,0)&(4,0) &(7,0)&(2,0)&(0,0)&(1,0)&(3,0)\\
(6,0)	&(6,0)&(7,0)&(5,0) &(4,0)&(0,0)&(1,0)&(3,0)&(2,0) \\
(7,0)	&(7,0)&(4,0)&(6,0)&(5,0)&(1,0)&(3,0)&(2,0)&(0,0)\\
(0,1)	&(0,1)& (2,1) &(1,1)& (3,1)&(6,1)&(5,1)&(4,1)&(7,1)\\
(1,1)	&(1,1)& (0,1)&(3,1)& (2,1)	&(5,1)&(4,1)&(7,1)	&(6,1)\\
(2,1)	&(2,1)&(3,1)&(0,1)	&(1,1)&(7,1)&(6,1)	&(5,1)&(4,1)\\
(3,1)	&(3,1)&(1,1)&(2,1)	&(0,1)&(4,1)&(7,1)&(6,1)&(5,1)\\
(4,1)	&(4,1)&(7,1)&(5,1)	&(6,1)&(0,1)&(2,1)&(3,1)&(1,1)\\
(5,1)	&(5,1)&(4,1)&(6,1)	&(7,1)&(1,1)&(0,1)&(2,1)&(3,1)\\
(6,1)	&(6,1)&(5,1)&(7,1)	&(4,1)&(3,1)&(1,1)&(0,1)&(2,1)\\
(7,1)	&(7,1)&(6,1)&(4,1)&(5,1)&(2,1)&(3,1)	&(1,1)&(0,1)\\
\end{tabular}

\begin{tabular}{c|cccccccc}
$\oplus$ & $(0,1)$ & $(1,1)$ & $(2,1)$ & $(3,1)$ & $(4,1)$ & $(5,1)$ & $(6,1) $ & $(7,1)$ \\ 
\hline \\
 & $(0,1)$ & $(1,1)$ & $(2,1)$ & $(3,1)$ & $(4,1)$ & $(5,1)$ 	& $(6,1)$ & $(7,1)$\\
& $(1,1)$ & $(3,1)$ & $(0,1)$ & $(2,1)$  & $(7,1)$  & $(4,1)$ & $(5,1)$ & $(6,1)$\\
 & $(2,1)$ &	$(0,1)$ & $(3,1)$ & $(1,1)$ & $(5,1)$ & $(6,1)$ & $(7,1)$ & $(4,1)$ \\
& (3,1)	&(2,1)&(1,1)&(0,1)&(6,1)&(7,1)	&(4,1)&(5,1)\\
& (4,1)&(5,1)&(7,1)&(6,1)&(3,1)	&(2,1)&(0,1)&(1,1)\\
& (5,1)&(6,1)&(4,1)&(7,1)&(2,1)&(0,1)	&(1,1)&(3,1)\\
& (6,1)&(7,1)&(5,1)&(4,1)&(0,1)	&(1,1)&(3,1)	&(2,1)\\
& (7,1)&(4,1)&(6,1)&(5,1)&(1,1)	&(3,1)&(2,1)&(0,1) \\
& (0,0)&(2,0)&(1,0) &(3,0)&(6,0)&(5,0)&(4,0)&(7,0) \\
& (1,0)&(0,0) &(3,0)&(2,0)	&(5,0)&(4,0)&(7,0)&(6,0)\\
& (2,0)&(3,0)&(0,0)&(1,0)&(7,0)&(6,0)&(5,0)&(4,0)\\
& (3,0)&(1,0)&(2,0)&(0,0)	&(4,0)&(7,0)&(6,0)&(5,0)\\
& (4,0)&(7,0)&(5,0)&(6,0)	&(0,0)&(2,0)&(3,0)&(1,0)\\
& (5,0)&(4,0)&(6,0)&(7,0)	&(1,0)&(0,0)&(2,0)&(3,0)\\
 & (6,0)&(5,0)&(7,0)&(4,0)&(3,0)	&(1,0)&(0,0)&(2,0)\\
 & (7,0)&(6,0)&(4,0)&(5,0)&(2,0)&(3,0)&(1,0)&(0,0)
\end{tabular}
}}
\caption{Cayley table for the gyro-group $ Dih(G_{8})$.} \label{tab: operation G8}
\end{center}
\end{table}

The set of vertices of the generalized dihedral gyro-group $ Dih(G_{8}) $ is 

as follows:
\begin{eqnarray*}
V_{(4, 0)} &= & \{ ((0,0), (3,0),(6,0),(4,0)), ((1,0), (2,0),(7,0), (5,0)), \\
 &  & ((3,1), (0,1), (6,1), (4,1)), ((1,1), (2,1), (5,1), (7,1))\}, \\
V_{(0, 1)} & = & \lbrace ((0,0), (0,1)),((1,1), (2,0)), ((3,1), (3,0)), ((4,0), (6,1)), ((0,5), (5,1)),\\
 & & ((4,1), (6, 0)), ((0,7), (7,1)), ((2,1), (1, 0)),\\
V_{(3, 0)} & = & \lbrace ((0,0), (3,0)),((1,0), (2,0)), ((4,0), (6,0)), ((7,0), (5,0)), ((3,1), (0,1)), \\
 &  & ((4,1), (6, 1)), ((5,1), (7,1)), ((2,1), (6, 1))\}
\end{eqnarray*}
For convenience, we display $ V_{(4, 0)} $ members with $ 0, 1, 2, 3$,  $ V_{(0, 1)} $ members with $ 4, 5, \cdots, 11 $, and $ V_{(3, 0)} $ members with $ 12, 13, \cdots, 19 $.

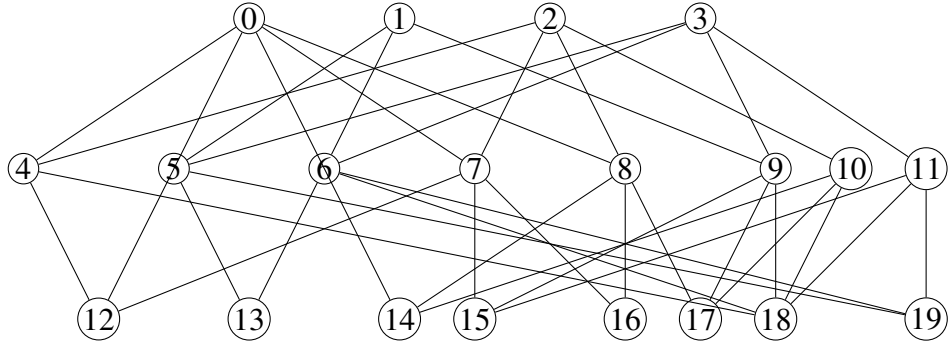
\begin{figure}[H]
\begin{center}
\begin{tikzpicture}
\node[circle,inner sep=.5pt,draw] (v0) at (0,0) {0};
\node[circle,inner sep=.5pt,draw] (v1) at (2,0) {1};
\node[circle,inner sep=.5pt,draw] (v2) at (4,0) {2};
\node[circle,inner sep=.5pt,draw] (v3) at (6,0) {3};
\node[circle,inner sep=.5pt,draw] (v4) at (-3,-2) {4};
\node[circle,inner sep=.5pt,draw] (v5) at (-1,-2) {5};
\node[circle,inner sep=.5pt,draw] (v6) at (1,-2) {6};
\node[circle,inner sep=.5pt,draw] (v7) at (3,-2) {7};
\node[circle,inner sep=.5pt,draw] (v8) at (5,-2) {8};
\node[circle,inner sep=.5pt,draw] (v9) at (7,-2) {9};
\node[circle,inner sep=.5pt,draw] (v10) at (8,-2) {10};
\node[circle,inner sep=.5pt,draw] (v11) at (9,-2) {11};
\node[circle,inner sep=.5pt,draw] (v12) at (-2,-4) {12};
\node[circle,inner sep=.5pt,draw] (v13) at (0,-4) {13};
\node[circle,inner sep=.5pt,draw] (v14) at (2,-4) {14};
\node[circle,inner sep=.5pt,draw] (v15) at (3,-4) {15};
\node[circle,inner sep=.5pt,draw] (v16) at (5,-4) {16};
\node[circle,inner sep=.5pt,draw] (v17) at (6,-4) {17};
\node[circle,inner sep=.5pt,draw] (v18) at (7,-4) {18};
\node[circle,inner sep=.5pt,draw] (v19) at (9,-4) {19};
\draw (v0) --(v4); 
\draw (v0) -- (v7);
\draw (v0) -- (v8);
\draw (v0) -- (v12);
\draw (v0) -- (v14);

\draw (v1) -- (v5);
\draw (v1) -- (v6);
\draw (v1) -- (v9);
\draw (v2) -- (v4);
\draw (v2) -- (v7);
\draw (v2) -- (v8);
\draw (v2) -- (v10);
\draw (v3) -- (v5);
\draw (v3) -- (v6);
\draw (v3) -- (v9);
\draw (v3) -- (v11);      
\draw (v4) -- (v12); 
\draw (v4) -- (v18);
\draw (v5) -- (v13);  
\draw (v5) -- (v19);
\draw (v6) -- (v13);  
\draw (v6) -- (v18);
\draw (v6) -- (v19);  
\draw (v7) -- (v12);  
\draw (v7) -- (v15); 
\draw (v7) -- (v16); 
\draw (v8) -- (v14); 
\draw (v8) -- (v16); 
\draw (v8) -- (v17);
 \draw (v9) -- (v15);
 \draw (v9) -- (v17);
 \draw (v9) -- (v18); 
 \draw (v10) -- (v14);
 \draw (v10) -- (v17);
 \draw (v10) -- (v18);
 \draw (v11) -- (v15);
 \draw (v11) -- (v18);
 \draw (v11) -- (v19);        
    \end{tikzpicture}
\caption{$\Phi(Dih(G_8),S=\{(3,0),(4,0),(0,1)\})$}
\end{center} 
\end{figure}
\end{exa}

\subsection{The Symmetry of Cayley graphs and $G$-graphs of gyro-groups}

One of the significant concept in graph theory is studying the automorphism group of a given graph. For instance in \cite{Ashrafi2019}, the automorphism groups of the Involution $G$-graph and Cayley graph are investigated. 
In this section we discuss about symmetry of Cayley graphs and $G$-graphs of gyro-groups. 
A graph  $ \Lambda $ is vertex-transitive if for all $ x, y \in V (\Lambda)$, there exists $f \in Aut(\Lambda)$ such that $f(x) = y$.
The edge-transitive graph is defined in the same way on the edges. A graph that is both a  vertex-transitive and an edge-transitive is called a symmetric graph.

In the previous sections, we talked about Cayley graph of gyro-groups and we saw that some features that exist in the Cayley graph of groups are not present in the Cayley graph of gyro-groups. In the following, we will examine another one of these features. 

We know that in Groups every Cayley graph is vertex-transitive, But this is not the case with gyro-groups.

\begin{exa}\cite{Bussaban} 
The Cayley graph of gyro-group $(G, \oplus) $ is defined in Example \ref{ma} with the generating set $ S =  \lbrace 1, 2, 3 \rbrace $ is not vertex-transitive.

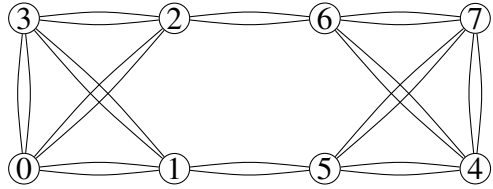
\begin{figure}[H]
\begin{center}
\begin{tikzpicture}
\node[circle,inner sep=.5pt,draw] (v1) at (0,0) {3};
\node[circle,inner sep=.5pt,draw] (v2) at (2,0) {2};
\node[circle,inner sep=.5pt,draw] (v3) at (4,0) {6};
\node[circle,inner sep=.5pt,draw] (v4) at (6,0) {7};
\node[circle,inner sep=.5pt,draw] (v5) at (6,-2) {4};
\node[circle,inner sep=.5pt,draw] (v6) at (4,-2) {5};
\node[circle,inner sep=.5pt,draw] (v7) at (2,-2) {1};
\node[circle,inner sep=.5pt,draw] (v8) at (0,-2) {0};
\draw (v1) edge[me=2] (v2); 
\draw (v2) edge[me=2] (v3); 
 \draw (v3) edge[me=2] (v4); 
 \draw (v4) edge[me=2] (v5); 
 \draw (v5) edge[me=2] (v6); 
\draw (v6) edge[me=2] (v7);  
\draw (v7) edge[me=2] (v8); 
\draw (v8) edge[me=2] (v1); 
 \draw (v1) edge[me=2] (v7); 
 \draw (v2) edge[me=2] (v8);  
  \draw (v3) edge[me=2] (v5); 
 \draw (v4) edge[me=2] (v6);
   \end{tikzpicture}
\caption{$Cay(G,\{1,2,3\})$}
\end{center} 
\end{figure}
\end{exa}

\begin{cor}
 We conclude that Cayley graph of a gyro-group need not be a vertex-transitive graph. In fact, it is still an open question " whether the Cayley graph of gyro-groups is transitive?"
\end{cor}

\begin{cor}
From the above discussion, we conclude that the Cayley graph of gyro-groups is not necessarily symmetric.
\end{cor}

In \cite{Maungchang, Maungchang2021} a sufficient condition is given for an $L$-Cayley graph and a $ R $-Cayley graph to be transitive.

\begin{thm}\cite{Maungchang}
Let $ (G, \oplus ) $ be a finite gyro-group and let $ S $ be a symmetric subset of $ G $. If $gyr[g, s]$ is the identity map for all $ g \in G$ and $s \in S  $  then $ L$-$Cay(G,S)$ is transitive.
\end{thm}

\begin{thm} \label{ma1}
Let $ (G, \oplus ) $ be a finite gyro-group and let $ S $ be a symmetric subset of $ G $. If $gyr [g, g^{\prime}] $ is the identity map for all $ g , g^{\prime} \in G$,  then $R$-$Cay(G,S)$ is transitive.
\end{thm}

In the following, we will examine the symmetry of $G$-graphs of gyro-graphs. 
Next, we have the following definition.
\begin{defn}
Let $ G $ be a gyro-group, $ S $ a subset of $ G $  so that $ S \neq \emptyset $. If $ S $ be left generating set, i.e $(S \rangle =G $, we call  the $G$-graph of $(G,S)  $ as

 $L$-$ G $-graph $ (G, S)$. and $S$ be right generating set, we say $R$-$G $-graph $(G,S)$.
\end{defn}

\begin{exa}
Consider the gyro-group $G_{16} $. let $ S = \lbrace 1,2,3 \rbrace $ be left generating set for $ G_{16} $. It is clear that S is symmetric.
Then we choose 3 arbitrary member $ 5, 6 \in G_{16} $ and $ 2 \in S  $. We have:
$$  gyr[5, 2](6)= \ominus (5 \oplus 2) \oplus (5 \oplus (2 \oplus 6) )= 6.$$
 On the other hand $ gyr[g,s]$ is the identity map for all $ g \in G$, $  s\in S $.  
 Then $ L $-$ G $-graph $  (G, S)$ is vertex-transitive.

\begin{table}[H]
\begin{center}
{\small{
\begin{tabular}{c|cccccccccccccccc}
 $\oplus$ & 0 & 1 & 2 & 3 & 4 & 5 & 6 & 7 & 8 & 9 & 10 & 11 & 12 & 13 & 14 & 15 \\ 
\hline 
0 & 0 & 1 & 2 & 3 & 4 & 5 & 6 & 7 & 8 & 9 & 10 & 11 & 12 & 13 & 14 & 15 \\
1 & 1 & 0 & 3 & 2 & 5 & 4 & 7 & 6 & 9 & 8 & 11 & 10 & 13 & 12 & 15 & 14 \\
2 & 2 & 3 & 1 & 0 & 6 & 7 & 5 & 4 & 11 & 10 & 8 & 9 & 15 & 14 & 12 & 13 \\
3 & 3 & 2 & 0 & 1 & 7 & 6 & 4 & 5 & 10 & 11 & 9 & 8 & 14 & 15 & 13 & 12 \\
4 & 4 & 5 & 6 & 7 & 3 & 2 & 0 & 1 & 15 & 14 & 12 & 13 & 9 & 8 & 11 & 10 \\
5 & 5 & 4 & 7 & 6 & 2 & 3 & 1 & 0 & 14 & 15 & 13 & 12 & 8 & 9 & 10 & 11 \\
6 & 6 & 7 & 5 & 4 & 0 & 1 & 2 & 3 & 13 & 12 & 15 & 14 & 10 & 11 & 9 & 8 \\
7 & 7 & 6 & 4 & 5 & 1 & 0 & 3 & 2 & 12 & 13 & 14 & 15 & 11 & 10 & 8 & 9 \\
8 & 8 & 9 & 10 & 11 & 12 & 13 & 14 & 15 & 0 & 1 & 2 & 3 & 4 & 5 & 6 & 7 \\
9 & 9 & 8 & 11 & 10 & 13 & 12 & 15 & 14 & 1 & 0 & 3 & 2 & 5 & 4 & 7 & 6 \\
10 & 10 & 11 & 9 & 8 & 14 & 15 & 13 & 12 & 3 & 2 & 0 & 1 & 7 & 6 & 4 & 5 \\
11 & 11 & 10 & 8 & 9 & 15 & 14 & 12 & 13 & 2 & 3 & 1 & 0 & 6 & 7 & 5 & 4 \\
12 & 12 & 13 & 14 & 15 & 11 & 10 & 8 & 9 & 6 & 7 & 5 & 4 & 0 & 1 & 2 & 3 \\
13 & 13 & 12 & 15 & 14 & 10 & 11 & 9 & 8 & 7 & 6 & 4 & 5 & 1 & 0 & 3 & 2 \\
14 & 14 & 15 & 13 & 12 & 8 & 9 & 10 & 11 & 4 & 5 & 6 & 7 & 3 & 2 & 0 & 1 \\
15 & 15 & 14 & 12 & 13 & 9 & 8 & 11 & 10 & 5 & 4 & 7 & 6 & 2 & 3 & 1& 0\\

\end{tabular}
\caption{The addition table of the gyro-group $G_{16}$.}
\label{Ta:G16addition}
}}
\end{center}
\end{table}
\end{exa}

\begin{thm} \label{m1}
Let $ G $ be a finite gyro-group and $ S $ be a symmetric subset of $ G $. If  $ gyr[g,s]$ is the identity map for all $ g \in G$, $  s\in S $. Then $ L $-$ G $-graph $  (G, S)$ is vertex-transitive.
\end{thm}
\begin{proof}
In fact, we prove that the condition on $ gyr[g,s] $ causes right addition by any $ g \in G$, be automorphisms on $ L $-$ G $-graph $  (G, S)$.
Let $ (u) \oplus x $ and $(v) \oplus y  $ be two vertices in $ L $-$ G $-graph $  (G, S)$, then
$$ ( \langle u \rangle \oplus x  \bigcap  \langle v \rangle \oplus x )= p \geq 1.$$
Now, Let $ \langle w \rangle \oplus x  $ and $ \langle z \rangle \oplus x  $ is adjacent in  $ L $-$ G $-graph $  (G, S)$. Therefor 
$$ ( \langle w \rangle \oplus x  \bigcap    \langle z \rangle \oplus y ) = m,$$
 where $m$ can be single point or a cyclic. We add $ g $ to both sides on the right. Then we have:
 \begin{eqnarray*}
 & & ( \langle w \rangle \oplus x  \bigcap    \langle z \rangle \oplus y ) \oplus g = \underbrace{m \oplus g} = l\\
& = &  ( ( \langle w \rangle \oplus x ) \oplus g )  \bigcap     ( (\langle z \rangle \oplus y )  \oplus g )= l \\
& =& ( (x, w_{1} \oplus x,\cdots ,  w_{n} \oplus x ) \oplus g ) \bigcap (( y,  z_{1} \oplus x,\cdots ,  z_{k} \oplus y ) \oplus g )= l  \\
& = & ((x \oplus g, (w_{1} \oplus x) \oplus g, \cdots, (w_{n} \oplus x) \oplus g )) \\
 & \bigcap & ( (y  \oplus g, (z_{1} \oplus y) \oplus g, \cdots, (z_{k} \oplus y) \oplus g ) ) = l \\
& = & ((x \oplus g, w_{1} \oplus (x \oplus gyr [w_{1}, x]), \cdots, w_{n} \oplus (x \oplus gyr [w_{n}, x]) = l \\
 & \bigcap &   ((y \oplus g, z_{1} \oplus (y \oplus gyr [z_{1}, y]), \cdots, z_{k} \oplus (y \oplus gyr [z_{k}, y]) = l 
 \end{eqnarray*}
since  $ gyr[g,s]$ is the identity map for all $ g \in G$, $  s\in S $ then:
\begin{eqnarray*}
 & = & ((x \oplus g, w_{1} \oplus  (x \oplus g ), \cdots, w_{n} \oplus (x \oplus g) ) \\
  & \bigcap &  ((y \oplus g, z_{1} \oplus  (y\oplus g ), \cdots, z_{k} \oplus (y \oplus g) ) =l \\
& = & ((w_{1}, \cdots, w_{n}) \oplus (x \oplus g)) \bigcap  ((z_{1}, \cdots, z_{k}) \oplus (y \oplus g)) = l  \\
 & = & ((w) \oplus (x \oplus g)) \bigcap ((z) \oplus (y \oplus g)) = l . 
\end{eqnarray*}
Then $ ((w) \oplus (x \oplus g)) $ and  $ ((z) \oplus (y \oplus g)) $ are adjacent. Hence permutation $ g_{s}= G \longmapsto  G $ where $ x \longrightarrow s \oplus x $
 is an automorphism that send $(u) \oplus x  $ to $ (v) \oplus y $. Therefore $ L$-$G$-graph is vertex-transitive.
\end{proof}

The same theorem holds for $R$-$G$- graphs.

\begin{thm}
Let $ G $ be a finite gyro-group with  a symmetric non-empty subset $S$. Assume that  $ gyr [g, g^{\prime}] (S)= S$ is the identity map for all $ g, g^{\prime} \in G$,  then $ R $-$G$-graph $(G, S)$ is vertex-transitive.
\end{thm}

\begin{proof}
Proof is similar to previous theorem.
\end{proof}

\begin{exa}
	Consider $ ( G_{16}, S= \lbrace 8, 9\rbrace )$, where  $ S $ is symmetric, $ (S \rangle =G $ and $gyr[g, g^{\prime} ](S)=S $ for $ g, g^{\prime} \in G $. We see $ R $-$ G $-graph $  (G, S)$ is vertex-transitive.
	
	\begin{figure}[H]
		\begin{center}
			{\footnotesize
				\begin{tikzpicture}
					\node (v1) at (-1,0) {$(0,8)$};
					\node (v2) at (1,0) {$(1,9)$};
					\node (v3) at (-1,-2) {$(0,9)$};
					\node (v4) at (1,-2) {$(1,8)$};
					\node (v5) at (2,0) {$(2,10)$};
					\node (v6) at (4,0) {$(3,11)$};
					\node (v7) at (2,-2) {$(2,11)$};
					\node (v8) at (4,-2) {$(3,10)$};
					\node (v9) at (5,0) {$(4,12)$};
					\node (v10) at (7,0) {$(5,13)$};
					\node (v11) at (5,-2) {$(4,13)$};
					\node (v12) at (7,-2) {$(5,12)$};
					\node (v13) at (8,0) {$(6,14)$};
					\node (v14) at (10,0) {$(7,15)$};
					\node (v15) at (8,-2) {$(6,15)$};
					\node (v16) at (10,-2) {$(7,14)$};
					\draw (v1) --(v3);
					\draw (v1) --(v4);
					\draw (v2)-- (v3);
					\draw (v2) -- (v4);
					\draw (v5) -- (v7);
					\draw (v5) -- (v8);
					\draw (v6) -- (v7);
					\draw (v6) -- (v8);
					\draw (v9) --(v11);
					\draw (v9) -- (v12);
					\draw (v10)-- (v11);
					\draw (v10) -- (v12);
					\draw (v13) --(v15);
					\draw (v13)-- (v16);
					\draw (v14) -- (v15);
					\draw (v14) -- (v16);
				\end{tikzpicture}
			}
			\caption{Cayley graph of $G_{16}$}
		\end{center}
	\end{figure}
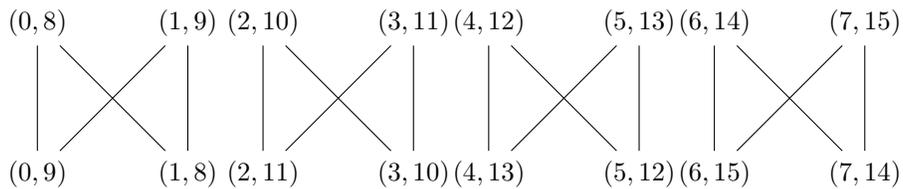
	
\end{exa}

\begin{rem}
The converse of the previous theorem is not true. For example:

 $( G_{8} , S=  \lbrace 1, 3 \rbrace )$ is vertex-transitive but some right additions are not the automorphism.
\end{rem}

In the following, we have the following definition, which will be necessary to prove the following theorems and propositions.

\begin{defn}\cite{Bretto-2007}
	The line graph $ L(\Gamma) $  associated with a simple graph $ \Gamma $ has for vertices the
	edges of 
	$ \Gamma $, two vertices being adjacent if and only if the corresponding edges in $ \Gamma $  are adjacent.
\end{defn}

\begin{prop}\label{mani}
Let $ G $ be a gyro-group and $ S = \lbrace s, t \rbrace $ with $ (S \rangle =G $ and

 $ \langle s \rangle \cap  \langle t \rangle =\lbrace 0 \rbrace $. Then $ \Gamma $= $ L $-$ G $-graph $  (G, S)$  is a simple graph and 

$ L(\Gamma) \simeq  Cay(G, A)$ where $ A= (\langle s \rangle \cup \langle t\rangle )\setminus \lbrace 0 \rbrace $.
\end{prop}

\begin{proof}
Let $ V $ the set of vertices of $  L(\Gamma) $. Then every member of the $ V $ is an edge of $ \Gamma $. i.e,
\begin{center}
$ a= ( [ (s) \oplus x, (t) \oplus y], u) $ where $ \langle s\rangle \oplus x  \cap  \langle t \rangle \oplus y= \lbrace u \rbrace $.
\end{center}

If $a=\theta (u) = \theta (u^{\prime})=a^{\prime}  $ we have:
\begin{center}
$([(s\rangle \oplus x, (t \rangle \oplus y ], u) = ([(s\rangle \oplus x^{\prime}, (t \rangle \oplus y ], u^{\prime})  $
\end{center}
As a result $ u=u^{\prime} $. Therefore $ \theta $ is one-to-one. And it is clear that $ \theta $ is surjective.

\end{proof}

\begin{exa}
Consider the gyro-group $G_8$ with the generating set $S=\{1,3\}$. The corresponding $G$-graph of this group is as follow.
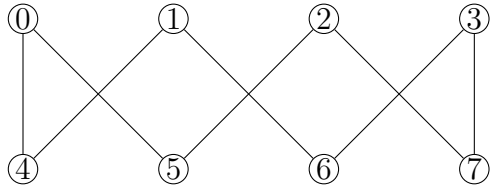
\begin{figure}[H]
\begin{center}
\begin{tikzpicture}
\node [circle,inner sep=.5pt,draw] (v1) at (0,0) {$0$};
\node [circle,inner sep=.5pt,draw] (v2) at (2,0) {$1$}; 
\node [circle,inner sep=.5pt,draw] (v3) at (4,0) {$ 2 $};
 \node [circle,inner sep=.5pt,draw] (v4) at (6,0) {$ 3 $}; 
 \node [circle,inner sep=.5pt,draw] (v5) at (6,-2) {$ 7 $};
\node [circle,inner sep=.5pt,draw] (v6) at (4,-2) {$6$}; 
\node [circle,inner sep=.5pt,draw] (v7) at (2,-2) {$5$};
 \node [circle,inner sep=.5pt,draw] (v8) at (0,-2) {$4$};
  \draw (v1) --(v8);  \draw (v1)-- (v7); \draw (v2) -- (v8); \draw (v2) --(v6); \draw (v3) -- (v5); \draw (v3) -- (v7); \draw(v4) -- (v5); \draw(v4)--(v6);
\end{tikzpicture}
\caption{$\Phi(G_8, \{1,3\})$.}
\end{center}
\end{figure}
\end{exa}

\begin{exa}
Consider $ G_{8}, S=\lbrace 1,3 \rbrace $, such that each edge $n-m$ is denoted by a letter as follows, then 
 $$ 04=a, 05=e,16=g, 27= p,14= f, 25= l, 36=q,  37=r ,$$ 
\begin{figure}[H]
\begin{center}
\begin{tikzpicture}
\node [circle,inner sep=.5pt,draw] (v1) at (0,0) {$a$};
\node [circle,inner sep=.5pt,draw] (v2) at (2,0) {$e$}; 
\node [circle,inner sep=.5pt,draw] (v3) at (4,0) {$ l $};
 \node [circle,inner sep=.5pt,draw] (v4) at (6,0) {$ p $}; 
 \node [circle,inner sep=.5pt,draw] (v5) at (6,-2) {$ r $};
\node [circle,inner sep=.5pt,draw] (v6) at (4,-2) {$q $}; 
\node [circle,inner sep=.5pt,draw] (v7) at (2,-2) {$g$};
 \node [circle,inner sep=.5pt,draw] (v8) at (0,-2) {$f$};
  \draw (v1) --(v2);  \draw (v2)-- (v3); \draw (v3) -- (v4); \draw (v4) --(v5); \draw (v5) -- (v6); \draw (v6) -- (v7); \draw(v7) -- (v8); \draw(v8)--(v1);
\end{tikzpicture}
\caption{$L(\Phi(G_8, \{1,3\}))$.}
\end{center}
\end{figure}
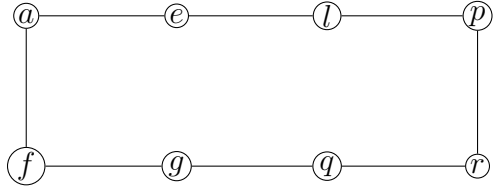

Finally, we drew the Cayley graph $( G_{8}, A= (\langle1 \rangle \cup  \langle 3 \rangle) \setminus \lbrace 0 \rbrace ) $  

\begin{figure}[H]
\begin{center}
\begin{tikzpicture}
\node [circle,inner sep=.5pt,draw] (v1) at (0,0) {$(4)$};
\node [circle,inner sep=.5pt,draw] (v2) at (2,0) {$(5)$}; 
\node [circle,inner sep=.5pt,draw] (v3) at (4,0) {$(1)$};
 \node [circle,inner sep=.5pt,draw] (v4) at (6,0) {$(0)$}; 
 \node [circle,inner sep=.5pt,draw] (v5) at (6,-2) {$(3)$};
\node [circle,inner sep=.5pt,draw] (v6) at (4,-2) {$(2)$}; 
\node [circle,inner sep=.5pt,draw] (v7) at (2,-2) {$(6)$};
 \node [circle,inner sep=.5pt,draw] (v8) at (0,-2) {$(7)$};
  \draw (v1) --(v2);  \draw (v2)-- (v3); \draw (v3) -- (v4); \draw (v4) --(v5); \draw (v5) -- (v6); \draw (v6) -- (v7); \draw(v7) -- (v8); \draw(v8)--(v1);
\end{tikzpicture}
\caption{$Cay(G_8, A=( \{ \langle 1 \rangle \cup \langle 3 \rangle ) \setminus \{0\})$.}
\end{center}
\end{figure}
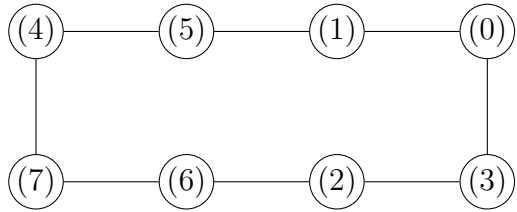

As we see, the line graph ($G$-graph  $  \left( G_{8}, \lbrace 1, 3 \rbrace  \right) $)  is isomorphic with Cayley graph  $ \left(  G_{8}, A=( \langle1 \rangle \cup  \langle 3 \rangle) \setminus \lbrace 0 \rbrace ) \right)  $.

\end{exa}

\begin{prop}\cite{Newman}  \label{ma2}
A connected graph is edge-transitive if and only if the line graph of its vertices is transitive.
\end{prop}

After the preliminaries, we will come to the most important point of this section.

\begin{thm}
Let $ G $ be a gyro-group and $ S = \lbrace s, t \rbrace $ with $ (S \rangle = G $ and

 $ \langle s \rangle \cap  \langle t \rangle = \lbrace 0 \rbrace $ and also $ gyr[s, g] $ is identity map for $ g \in G $ and $ s \in S $, then $ L $-$ G $-graph $  (G, S)$ is symmetric.
\end{thm}

\begin{proof}
According to theorem \ref{m1}, this graph is vertex-transitive. According to proposition \ref{mani} line graph of this graph is isomorphic with Cayley graph and is vertex-transitive. By using of the  proposition \ref{ma2}, we can see this graph is edge-transitive. This completes the proof.
\end{proof}

\newpage

\noindent{\bf Acknowledgement.} 
 The research of the second author is partially supported by the University of Kashan under grant no 1311841.

\end{document}